\documentclass[11pt]{article}

\usepackage{latexsym}
\usepackage{amsmath}
\usepackage{amsfonts}
\usepackage{graphicx,epsfig}

\newtheorem{Thm}{Theorem}[section]

\newtheorem{lemma}[Thm]{Lemma}

\newtheorem{proposition}[Thm]{Proposition}
\newtheorem{definition}[Thm]{Definition}

\newtheorem{theorem}[Thm]{Theorem}

\newcommand{\bitem}{\begin{itemize}}
\newcommand{\eitem}{\end{itemize}}
\newcommand{\benum}{\begin{enumerate}}
\newcommand{\eenum}{\end{enumerate}}
\newcommand{\beq}{\begin{equation}}
\newcommand{\eeq}{\end{equation}}
\newcommand{\ip}[2]{\langle#1,#2\rangle}

\newcommand{\absip}[2]{| \langle#1,#2\rangle |}
\newcommand{\norm}[1]{\|#1\|}

\def\cC{{\mathcal{C}}}

\def\cH{{\mathcal{H}}}

\newcommand{\cS}{{\mathcal S}}
\newcommand{\cN}{{\mathcal N}}
\newcommand{\cP}{{\mathcal P}}

\newcommand{\qed}{$\Box$}

\usepackage[caption=false]{subfig}
\usepackage{url}

\def\cH{\mathcal{H}}

\begin{document}

\title{Data Separation by Sparse Representations}

\author{Gitta Kutyniok\\[1ex]Institute of Mathematics, University of Osnabr\"uck, 49069 Osnabr\"uck, Germany}


\maketitle

\begin{abstract}
Recently, sparsity has become a key concept in various areas of applied mathematics, computer science,
and electrical engineering. One application of this novel methodology is the separation of data, which
is composed of two (or more) morphologically distinct constituents. The key idea is to carefully select
representation systems each providing sparse approximations of one of the components.
Then the sparsest coefficient vector representing the data within the composed -- and therefore highly
redundant -- representation system is computed by $\ell_1$ minimization or thresholding. This automatically
enforces separation.

This paper shall serve as an introduction to and a survey about this exciting area of research as well as a
reference for the state-of-the-art of this research field.
\end{abstract}



\vspace{.1in}
\vspace{.1in}

\vspace{.1in}
{\bf Acknowledgements.}
The author would like to thank Ronald Coifman, Michael Elad, and
Remi Gribonval for various discussions on related topics, and Wang-Q
Lim for producing Figures \ref{fig:circle}, \ref{fig:neuron}, and
\ref{fig:barbara}. Special thanks go to David Donoho for a great
collaboration on topics in this area and enlightening debates, and
to Michael Elad for very useful comments on an earlier version of
this survey. The author is also grateful to the Department of
Statistics at Stanford University and the Mathematics Department at
Yale University for their hospitality and support during her visits.
She acknowledges support by Deutsche Forschungsgemeinschaft (DFG)
Heisenberg fellowship KU 1446/8, DFG Grant SPP-1324 KU 1446/13, and
DFG Grant KU 1446/14.
\vspace{.1in}

\pagebreak

\section{Introduction}

Over the last years, scientists face an ever growing deluge of data,
which needs to be transmitted, analyzed, and stored. A close
analysis reveals that most of these data might be classified as
multimodal data, i.e., being composed of distinct subcomponents.
Prominent examples are audio data, which might consist of a
superposition of the sounds of different instruments, or imaging
data from neurobiology, which is typically a composition of the soma
of a neuron, its dendrites, and its spines. In both these exemplary
situations, the data has to be separated into appropriate single
components for further analysis. In the first case, separating the
audio signal into the signals of the different instruments is a
first step to enable the audio technician to obtain a musical score
from a recording. In the second case, the neurobiologist might aim
to analyze the structure of dendrites and spines separately for the
study of Alzheimer specific characteristics. Thus data separation is
often a crucial step in the analysis of data.

As a scientist, three fundamental problems immediately come to one's
mind:
\begin{enumerate}
\item[(P1)] What is a mathematically precise meaning of the vague term `distinct components'?
\item[(P2)] How do we separate data algorithmically?
\item[(P3)] When is separation possible at all?
\end{enumerate}
To answer those questions, we need to first understand the key
problem in data separation. In a very simplistic view, the essence
of the problem is as follows: Given a composed signal $x$ of the
form $x = x_1 + x_2$, we aim to extract the unknown components $x_1$
and $x_2$ from it. Having one known data and two unknowns obviously
makes this problem underdetermined. Thus, the novel paradigm of
sparsity -- appropriately utilized -- seems a perfect fit for
attacking data separation, and this chapter shall serve as both an
introduction into this intriguing application of sparse
representations as well as a reference for the state-of-the-art of
this research area.

\subsection{Morphological Component Analysis}

Intriguingly, when considering the history of Compressed Sensing,
the first mathematically precise result on recovery of sparse
vectors by $\ell_1$ minimization is related to a data separation
problem: The separation of sinusoids and spikes in \cite{DS89,DH01}.
Thus it might be considered a milestone in the development of
Compressed Sensing. In addition, it reveals a surprising connection
with uncertainty principles.

The general idea allowing separation in \cite{DS89,DH01} was to
choose two bases or frames $\Phi_1$ and $\Phi_2$ adapted to the two
components to be separated in such a way that $\Phi_1$ and $\Phi_2$
provide a sparse representation for $x_1$ and $x_2$, respectively.
Searching for the sparsest representation of the signal in the
combined (highly overcomplete) dictionary $[\:\Phi_1\:|\:\Phi_2\:]$
should then intuitively enforce separation provided that $x_1$ does
not have a sparse representation in $\Phi_2$ and that $x_2$ does not
have a sparse representation in $\Phi_1$. This general concept was
later -- in the context of image separation, but the term seems to
be fitting in general -- coined {\em Morphological Component
Analysis} \cite{SMBED05}.

This viewpoint now measures the morphological difference between
components in terms of the incoherence of suitable sparsifying bases
or frames $\Phi_i$, thereby giving one possible answer to (P1); see
also the respective chapters in the novel book \cite{Sta10}. One possibility
for measuring incoherence is the {\em mutual coherence}. We will
however see in the sequel that there exist even more appropriate
coherence notions, which provide a much more refined measurement of
incoherence specifically adapted to measuring morphological
difference.

\subsection{Separation Algorithms}

Going again back in time, we observe that far before \cite{DH01},
Coifman, Wickerhauser, and co-workers already presented very
inspiring empirical results on the separation of image components
using the idea of Morphological Component Analysis, see \cite{CW93}.
After this, several techniques to actually compute the sparsest
expansion in a composed dictionary $[\:\Phi_1\:|\:\Phi_2\:]$ were
introduced. In \cite{MZ93}, Mallat and Zhang developed {\em Matching
Pursuit} as one possible methodology. The study by Chen, Donoho, and
Saunders in \cite{CDS98} then revealed that the $\ell_1$ norm has a
tendency to find sparse solutions when they exist, and coined this
method {\em Basis Pursuit}.

As explained before, data separation by Morphological Component
Analysis -- when suitably applied -- can be reduced to a sparse
recovery problem. To solve this problem, there nowadays already
exist a variety of utilizable algorithmic approaches; thereby
providing a general answer to (P2). Such approaches include, for
instance, a canon of greedy-type algorithms. Most of the theoretical separation results however consider
$\ell_1$ minimization as the main separation technique, which is
what we will also mainly focus on in this chapter.

\subsection{Separation Results}

As already mentioned, the first mathematically precise result was
derived in \cite{DH01} and solved the problem of separation of
sinusoids and spikes. After this `birth of sparse data separation',
a deluge of very exciting results started. One direction of research
are general results on sparse recovery and Compressed Sensing; here
we would like to cite the excellent survey paper \cite{BDE09}.

Another direction continued the idea of sparse data separation
initiated in \cite{DH01}. In this realm, the most significant
theoretical results might be considered firstly the series of papers
\cite{EB02,DE03}, in which the initial results from \cite{DH01} are
extended to general composed dictionaries, secondly the paper
\cite{GN03}, which also extends results from \cite{DH01} though with
a different perspective, and thirdly the papers \cite{BGN08} and
\cite{DK10}, which explore the clustering of the sparse coefficients
and the morphological difference of the components encoded in it.

We also wish to mention the abundance of empirical work showing that
utilizing the idea of sparse data separation often gives very
compelling results in practice, as examples, we refer to the series
of papers on applications to astronomical data
\cite{BSFMD07,SMBED05,SED04}, to general imaging data
\cite{MAC02,ESQD05,SED05}, and to audio data \cite{GB03,KT09}.

Let us remark that also the classical problem of denoising can be
regarded as a separation problem, since we aim to separate a signal
from noise by utilizing the characteristics of the signal family and
the noise. However, as opposed to the separation problems discussed
in this chapter, denoising is not a `symmetric' separation task,
since the characterization of the signal and the noise are very
different.

\subsection{Design of Sparse Dictionaries}
\label{subsec:design}

For satisfactorily answering (P3), one must also raise the question
of how to find suitable sparsifying bases or frames for given
components. This search for `good' systems in the sense of sparse
dictionaries can be attacked in two ways, either non-adaptively or
adaptively.

The first path explores the structure of the component one would
like to extract, for instance, it could be periodic such as
sinusoids or anisotropic such as edges in images. This typically
allows one to find a suitable system among the already very well
explored representation systems such as the Fourier basis, wavelets,
or shearlets, to name a few. The advantage of this approach is the
already explored structure of the system, which can hence be
exploited for deriving theoretical results on the accuracy of
separation, and the speed of associated transforms.

The second path uses a training set of data similar to the
to-be-extracted component, and `learns' a system which best
sparsifies this data set. Using this approach customarily referred
to as {\em dictionary learning}, we obtain a system extremely well
adapted to the data at hand; as the state-of-the-art we would like
to mention the K-SVD algorithm introduced by Aahron, Elad, and
Bruckstein in \cite{AEB06}; see also \cite{DS09} for a `Compressed
Sensing' perspective to K-SVD. Another appealing dictionary training
algorithm, which should be cited is the method of optimal directions
(MOD) by Engan et al. \cite{EAHH99}. The downside however is the
lack of a mathematically exploitable structure, which makes a
theoretical analysis of the accuracy of separation using such a
system very hard.

\subsection{Outline}

In Section \ref{sec:estimates}, we discuss the formal mathematical
setting of the problem, present the nowadays already considered
classical separation results, and then discuss more recent results
exploiting the clustering of significant coefficients in the
expansions of the components as a means to measure their
morphological difference. We conclude this section by revealing a
close link of data separation to uncertainty principles. Section
\ref{sec:signals} is then devoted to both theoretical results as
well as applications for separation of 1D signals, elaborating, in
particular, on the separation of sinusoids and spikes. Finally,
Section \ref{sec:images} focuses on diverse questions concerning
separation of 2D signals, i.e., images, such as the separation of
point- and curvelike objects, again presenting both application
aspects as well as theoretical results.


\section{Separation Estimates}
\label{sec:estimates}

As already mentioned in the introduction, data separation can be
regarded within the framework of underdetermined problems. In this
section, we make this link mathematically precise. Then we discuss
general estimates on the separability of composed data, firstly
without any knowledge of the geometric structure of sparsity
patterns, and secondly, by taking known geometric information into
account. A revelation of the close relation with uncertainty
principles concludes the section.

In Sections \ref{sec:signals} and \ref{sec:images}, we will then see
the presented general results and uncertainty principles in action,
i.e., applied to real-world separation problems.


\subsection{Relation with Underdetermined Problems}

Let $x$ be our signal of interest, which we for now consider as
belonging to some Hilbert space $\cH$, and assume that
\[
x = x_1^0 + x_2^0.
\]
Certainly, real data is typically composed of multiple components,
hence not only the situation of two components, but three or more is
of interest. We will however focus on the two-component situation to
clarify the fundamental principles behind the success of separating
those by sparsity methodologies. It should be mentioned though that,
in fact, most of the presented theoretical results can be extended
to the multiple component situation in a more or less
straightforward manner.

To extract the two components from $x$, we need to assume that --
although we are not given $x_1^0$ and $x_2^0$ -- certain
`characteristics' of those components are known to us. Such
`characteristics' might be, for instance, the pointlike structure of
stars and the curvelike structure of filaments in astronomical
imaging. This knowledge now enables us to choose two representation
systems, $\Phi_1$ and $\Phi_2$, say, which allow sparse expansions
of $x_1^0$ and $x_2^0$, respectively. Such representation systems
might be chosen from the collection of well-known systems such as
wavelets. A different possibility is to choose adaptively the
systems via dictionary learning procedures. This approach however
requires training data sets for the two components $x_1^0$ and
$x_2^0$ as discussed in Subsection \ref{subsec:design}.

Given now two such representation systems $\Phi_1$ and $\Phi_2$, we
can write $x$ as
\[
x = x_1^0 + x_2^0 = \Phi_1 c_1^0 + \Phi_2 c_2^0 =
[\:\Phi_1\:|\:\Phi_2\:] \left[ \begin{array}{c} c_1^0\\ c_2^0
\end{array}\right]
\]
with $\|c_1^0\|_0$ and $\|c_2^0\|_0$ `sufficiently small'. Thus, the
data separation problem has been reduced to solving the
underdetermined linear system \beq \label{eq:underdetermined} x =
[\:\Phi_1\:|\:\Phi_2\:] \left[ \begin{array}{c} c_1\\ c_2
\end{array}\right] \eeq for $[c_1, c_2]^T$. Unique recovery of the
original vector $[c_1^0, c_2^0]^T$ automatically extracts the
correct two components $x_1^0$ and $x_2^0$ from $x$, since
\[
x_1^0 = \Phi_1 c_1^0 \quad \mbox{and} \quad x_2^0 = \Phi_2 c_2^0.
\]
Ideally, one might want to solve \beq \label{eq:mingeneral0}
\min_{c_1,c_2} \norm{c_1}_0 + \norm{c_2}_0 \quad \textrm{s.t.} \quad
x = [\:\Phi_1\:|\:\Phi_2\:] \left[ \begin{array}{c} c_1\\ c_2
\end{array}\right], \eeq which however is an NP-hard problem. Instead one aims
to solve the $\ell_1$ minimization problem \beq
\label{eq:mingeneral} (\mbox{Sep}_s) \qquad  \min_{c_1,c_2}
\norm{c_1}_1 + \norm{c_2}_1 \quad \textrm{s.t.} \quad  x =
[\:\Phi_1\:|\:\Phi_2\:] \left[ \begin{array}{c} c_1\\ c_2
\end{array}\right]. \eeq The lower case `s' in $\mbox{Sep}_s$
indicates that the $\ell_1$ norm is placed on the synthesis side.
Other choices for separation are, for instance, greedy-type
algorithms. In this chapter we will
focus on $\ell_1$ minimization as the separation technique,
consistent with most known separation results from the literature.

Before discussing conditions on $[c_1^0, c_2^0]^T$ and
$[\:\Phi_1\:|\:\Phi_2\:]$, which guarantee unique solvability of
\eqref{eq:underdetermined}, let us for a moment debate whether
uniqueness is necessary at all. If $\Phi_1$ and $\Phi_2$ form bases,
it is certainly essential to recover $[c_1^0, c_2^0]^T$ uniquely
from \eqref{eq:underdetermined}. However, some well-known
representation systems are in fact redundant and typically
constitute Parseval frames such as curvelets or shearlets. Also,
systems generated by dictionary learning are normally highly
redundant. In this situation, for each possible separation \beq
\label{eq:exampleseparation} x = x_1 + x_2, \eeq there exist
infinitely many coefficient sequences $[c_1, c_2]^T$ satisfying \beq
\label{eq:redundancyproblem} x_1 = \Phi_1 c_1 \quad \mbox{and} \quad
x_2 = \Phi_2 c_2. \eeq Since we are {\em only} interested in the
correct separation and {\em not} in computing the sparsest
expansion, we can circumvent presumably arising numerical
instabilities when solving the minimization problem
\eqref{eq:mingeneral} by selecting a particular coefficient sequence
for each separation. Assuming $\Phi_1$ and $\Phi_2$ are Parseval
frames, we can exploit this structure and rewrite
\eqref{eq:redundancyproblem} as
\[
x_1 = \Phi_1 (\Phi_1^T x_1) \quad \mbox{and} \quad x_2 = \Phi_2
(\Phi_2^T x_2).
\]
Thus, for each separation \eqref{eq:exampleseparation}, we choose a
{\em specific} coefficient sequence when expanding the components in
the Parseval frames, in fact, we choose the {\em analysis sequence}.
This leads to the following different $\ell_1$ minimization problem
in which the $\ell_1$ norm is placed on the {\em analysis} rather
than the {\em synthesis} side: \beq \label{eq:newmingeneral}
(\mbox{Sep}_{a}) \qquad   \min_{x_1, x_2} \norm{\Phi_1^T x_1}_1 +
\norm{\Phi_2^T x_2}_1 \quad \textrm{s.t.} \quad  x = x_1 + x_2. \eeq

This new minimization problem can be also regarded as a mixed
$\ell_1$-$\ell_2$ problem, since the analysis coefficient sequence
is exactly the coefficient sequence which is minimal in the $\ell_2$
norm.


\subsection{General Separation Estimates}
\label{subsec:generalseparation}

Let us now discuss the main results of successful data separation,
i.e., stating conditions on $[c_1^0, c_2^0]^T$ and
$[\:\Phi_1\:|\:\Phi_2\:]$ for extracting $x_1^0$ and $x_2^0$ from
$x$. The strongest known general result was derived in 2003 by
Donoho and Elad \cite{DE03} and used the notion of mutual coherence.
Recall that, for a normalized frame $\Phi =(\varphi_{i})_{i \in I}$,
the {\em mutual coherence} of $\Phi$ is defined by
\[
\mu(\Phi) = \max_{i,j \in I, i \neq j}  | \langle \varphi_{i},
\varphi_{j} \rangle|.
\]
The result states the following.

\begin{theorem} [\cite{DE03}] \label{theo:generalbound}
Let $\Phi_1$ and $\Phi_2$ be two frames for a Hilbert space $\cH$,
and let $x \in \cH$, $x \neq 0$. If $x = [\Phi_1|\Phi_2]c$ and
\[
\|c\|_0 < \frac{1}{2} \left(1+\frac{1}{\mu([\Phi_1|\Phi_2])}\right),
\]
then the solution of the $\ell_1$ minimization problem
$(\mbox{Sep}_{s})$ stated in \eqref{eq:mingeneral} coincides with
the solution of the $\ell_0$ minimization problem stated in
\eqref{eq:mingeneral0}.
\end{theorem}

Before presenting the proof, we require some prerequisites. Firstly,
we need to introduce the so-called nullspace property.

\begin{definition}
Let $\Phi = (\varphi_i)_{i \in I}$ be a frame for a Hilbert space
$\cH$, and let $\cN(\Phi)$ denote the null space of $\Phi$. Then
$\Phi$ is said to have the {\em null space property of order $k$} if
\[
\|1_\Lambda d\|_1 < \frac12 \|d\|_1
\]
for all $d \in \cN(\Phi) \setminus \{0\}$ and for all sets $\Lambda
\subseteq I$ with $|\Lambda| \le k$.
\end{definition}

This notion provides a very useful characterization of the existence
of unique sparse solutions of the $\ell_1$ minimization problem
$(\mbox{Sep}_{s})$ stated in \eqref{eq:mingeneral}.

\begin{lemma} \label{lem:generalbound1}
Let $\Phi = (\varphi_i)_{i \in I}$ be a frame for a Hilbert space
$\cH$, and let $x \in \cH$. Then the following conditions are
equivalent.
\begin{enumerate}
\item[(i)] All vectors $c$ with $\|c\|_0 \le k$ are unique solutions of the $\ell_1$ minimization problem
$(\mbox{Sep}_{s})$ stated in \eqref{eq:mingeneral} (with $\Phi$
instead of $[\Phi_1|\Phi_2]$).
\item[(ii)] $\Phi$ satisfies the null space property of order $k$.
\end{enumerate}
\end{lemma}

\noindent {\bf Proof.}
 First, assume that (i) holds. Let $d \in
\cN(\Phi) \setminus \{0\}$ and $\Lambda \subseteq I$ with $|\Lambda|
\le k$ be arbitrary. Then, by (i), the sparse vector $1_\Lambda d$
is the unique minimizer of $\|c\|_1$ subject to $\Phi c = \Phi
(1_\Lambda d)$. Further, since $d \in \cN(\Phi) \setminus \{0\}$,
\[
\Phi(-1_{\Lambda^c} d) = \Phi(1_{\Lambda} d).
\]
Hence
\[
\|1_\Lambda d\|_1 < \|1_{\Lambda^c} d\|_1,
\]
or, in other words,
\[
\|1_\Lambda d\|_1 < \frac12 \|d\|_1,
\]
which implies (ii), since $d$ and $\Lambda$ were chosen arbitrarily.

Secondly, assume that (ii) holds, and let $c_1$ be a vector with
$\|c_1\|_0 \le k$ and support denoted by $\Lambda$. Further, let
$c_2$ be an arbitrary solution of $x = \Phi c$, and set
\[
d = c_2-c_1.
\]
Then
\[
\|c_2\|_1 - \|c_1\|_1 = \|1_{\Lambda^c} c_2\|_1 + \|1_{\Lambda}
c_2\|_1 - \|1_{\Lambda} c_1\|_1  \ge \|1_{\Lambda^c} d\|_1 -
\|1_\Lambda d\|_1.
\]
This term is greater than zero for any $d \neq 0$ if
\[
\|1_{\Lambda^c} d\|_1 > \|1_\Lambda d\|_1,
\]
or
\[
\frac12 \|d\|_1 > \|1_\Lambda d\|_1.
\]
This is ensured by (ii). Hence $\|c_2\|_1 > \|c_1\|_1$, and thus
$c_1$ is the unique solution of $(\mbox{Sep}_{s})$. This implies
(i).
\qed

Using this result, we next prove that a solution satisfying $\|c\|_0
< \frac{1}{2} \left(1+\frac{1}{\mu(\Phi)}\right)$ is the unique
solution of the $\ell_1$ minimization problem $(\mbox{Sep}_{s})$.

\begin{lemma} \label{lem:generalbound2}
Let $\Phi = (\varphi_i)_{i \in I}$ be a frame for a Hilbert space
$\cH$, and let $x \in \cH$. If $c$ is a solution of the $\ell_1$
minimization problem $(\mbox{Sep}_{s})$ stated in
\eqref{eq:mingeneral} (with $\Phi$ instead of $[\Phi_1|\Phi_2]$) and
satisfies
\[
\|c\|_0 < \frac{1}{2} \left(1+\frac{1}{\mu(\Phi)}\right),
\]
then it is the unique solution.
\end{lemma}

\noindent {\bf Proof.}
Let $d \in \cN(\Phi) \setminus \{0\}$, hence, in particular,
\[
\Phi d = 0;
\]
thus also \beq \label{eq:generalbound2_1} \Phi^\star \Phi d = 0.
\eeq Without loss of generality, we now assume that the vectors in
$\Phi$ are normalized.  Then, \eqref{eq:generalbound2_1} implies
that, for all $i \in I$,
\[
d_i = - \sum_{j \neq i} \ip{\varphi_i}{\varphi_j} d_j.
\]
Using the definition of mutual coherence $\mu(\Phi)$ (cf.
Subsection \ref{subsec:generalseparation}), we
obtain
\[
|d_i| \le \sum_{j \neq i} |\ip{\varphi_i}{\varphi_j}| \cdot |d_j|
\le \mu(\Phi) (\|d\|_1 - |d_i|),
\]
and hence
\[
|d_i| \le \left(1+\frac{1}{\mu(\Phi)}\right)^{-1} \|d\|_1.
\]
Thus, by the hypothesis on $\|c\|_0$ and for any $\Lambda \subseteq
I$ with $|\Lambda| = \|c\|_0$, we have
\[
\|1_\Lambda d\|_1 \le |\Lambda| \cdot
\left(1+\frac{1}{\mu(\Phi)}\right)^{-1} \|d\|_1 = \|c\|_0 \cdot
\left(1+\frac{1}{\mu(\Phi)}\right)^{-1} \|d\|_1 < \frac12 \|d\|_1.
\]
This shows that $\Phi$ satisfies the null space property of order
$\|c\|_0$, which, by Lemma \ref{lem:generalbound1}, implies that $c$
is the unique solution of $(\mbox{Sep}_{s})$.
\qed

We further prove that a solution satisfying $\|c\|_0 < \frac{1}{2}
\left(1+\frac{1}{\mu(\Phi)}\right)$ is also the unique solution of
the $\ell_0$-minimization problem.

\begin{lemma} \label{lem:generalbound3}
Let $\Phi = (\varphi_i)_{i \in I}$ be a frame for a Hilbert space
$\cH$, and let $x \in \cH$. If $c$ is a solution of the $\ell_0$
minimization problem stated in \eqref{eq:mingeneral0} (with $\Phi$
instead of $[\Phi_1|\Phi_2]$) and satisfies
\[
\|c\|_0 < \frac{1}{2} \left(1+\frac{1}{\mu(\Phi)}\right),
\]
then it is the unique solution.
\end{lemma}

\noindent {\bf Proof.}
By Lemma \ref{lem:generalbound2}, the hypotheses imply that $c$ is
the unique solution of the $\ell_1$ minimization problem
$(\mbox{Sep}_{s})$. Now, towards a contradiction, assume that there
exists some $\tilde{c}$ satisfying $x = \Phi \tilde{c}$ with
$\|\tilde{c}\|_0 \le \|c\|_0$. Then $\tilde{c}$ must satisfy
\[
\|\tilde{c}\|_0 < \frac{1}{2} \left(1+\frac{1}{\mu(\Phi)}\right).
\]
Again, by Lemma \ref{lem:generalbound2}, $\tilde{c}$ is the unique
solution of the $\ell_1$ minimization problem $(\mbox{Sep}_{s})$, a
contradiction.
\qed

These lemmata now immediately imply Theorem \ref{theo:generalbound}.\\[-0.3cm]

\noindent {\bf Proof [Proof of Theorem \ref{theo:generalbound}].}
Theorem \ref{theo:generalbound} follows from Lemmata
\ref{lem:generalbound2} and \ref{lem:generalbound3}.
\qed\\[-0.3cm]

Interestingly, in the situation of $\Phi_1$ and $\Phi_2$ being two
orthonormal bases the bound can be slightly strengthened. For the
proof of this result, we refer the reader to \cite{EB02}.

\begin{theorem} [\cite{EB02}] \label{theo:specialbound}
Let $\Phi_1$ and $\Phi_2$ be two orthonormal bases for a Hilbert
space $\cH$, and let $x \in \cH$. If $x = [\Phi_1|\Phi_2]c$ and
\[
\|c\|_0 < \frac{\sqrt{2}-0.5}{\mu([\Phi_1|\Phi_2])},
\]
then the solution of the $\ell_1$ minimization problem
$(\mbox{Sep}_{s})$ stated in \eqref{eq:mingeneral} coincides with
the solution of the $\ell_0$ minimization problem stated in
\eqref{eq:mingeneral0}.
\end{theorem}

This shows that in the special situation of two orthonormal bases,
the bound is nearly a factor of $2$ stronger than in the general
situation of Theorem \ref{theo:generalbound}.


\subsection{Clustered Sparsity as a Novel Viewpoint}
\label{subsec:clusteredsparsity}

In a concrete situation, we often have more information on the
geometry of the to-be-separated components $x_1^0$ and $x_2^0$. This
information is typically encoded in a particular clustering of the
non-zero coefficients if a suitable basis or frame for the expansion
of $x_1^0$ or $x_2^0$ is chosen. Think, for instance, of the tree
clustering of wavelet coefficients of a point singularity. Thus, it
seems conceivable that the morphological difference is encoded not
only in the incoherence of the two chosen bases or frames adapted to
$x_1^0$ and $x_2^0$, but in the interaction of the elements of those
bases or frames associated with the clusters of significant
coefficients. This should intuitively allow for weaker necessary
conditions for separation.

One possibility for a notion capturing this idea is the so-called
{\em joint concentration} which was introduced in \cite{DK10} with
concepts going back to \cite{DS89}, and was in between again revived
in \cite{DH01}. To provide some intuition for this notion, let
$\Lambda_1$ and $\Lambda_2$ be subsets of indexing sets of two
Parseval frames. Then the joint concentration measures the maximal
fraction of the total $\ell_1$ norm which can be concentrated on the
index set $\Lambda_1 \cup \Lambda_2$ of the combined dictionary.

\begin{definition}
Let $\Phi_1 =(\varphi_{1i})_{i \in I}$ and $\Phi_2
=(\varphi_{2j})_{j \in J}$ be two Parseval frames for a Hilbert
space $\cH$. Further, let $\Lambda_1 \subseteq I$ and $\Lambda_2
\subseteq J$. Then the {\em joint concentration} $\kappa= \hspace*{-0.1cm}
\kappa(\Lambda_1, \Phi_1; \Lambda_2, \Phi_2)$ is defined by
\[
\kappa(\Lambda_1, \Phi_1; \Lambda_2, \Phi_2) = \sup_x
\frac{\norm{1_{\Lambda_1} \Phi_1^T x}_1 + \norm{1_{\Lambda_2}
\Phi_2^T x}_1}{\norm{\Phi_1^T x}_1 + \norm{\Phi_2^T x}_1}.
\]
\end{definition}

One might ask how the notion of joint concentration relates to the
widely exploited, and for the previous result utilized mutual
coherence. For this, we first briefly discuss some derivations of
mutual coherence. A first variant better adapted to clustering of
coefficients was the {\em Babel function} first introduced in
\cite{DE03} and later in \cite{Tro04} under the label {\em
cumulative coherence function}, which, for a normalized frame $\Phi
=(\varphi_{i})_{i \in I}$ and some $m \in \{1, \ldots, |I|\}$ is
defined by
\[
\mu_B(m, \Phi) = \max_{\Lambda \subset I, |\Lambda|=m} \max_{j
\not\in \Lambda}  \sum_{i \in I} | \langle \varphi_{i}, \varphi_{j}
\rangle|.
\]
This notion was later refined in \cite{BGN08} by considering the
so-called {\em structured $p$-Babel function}, defined for some
family $\cS$ of subsets of $I$ and some $1 \le p < \infty$ by
\[
\mu_{sB}(\cS, \Phi) = \max_{\Lambda \in \cS} \left(\max_{j \not\in
\Lambda}  \sum_{i \in I} | \langle \varphi_{i}, \varphi_{j}
\rangle|^p\right)^{1/p}.
\]
Another variant, better adapted to data separation, is the {\em
cluster coherence} introduced in \cite{DK10}, whose definition we
now formally state. Notice that we do not assume that the vectors
are normalized.

\begin{definition} \label{def:cluster}
Let $\Phi_1 =(\varphi_{1i})_{i \in I}$ and $\Phi_2
=(\varphi_{2j})_{j \in J}$ be two Parseval frames for a Hilbert
space $\cH$, let $\Lambda_1 \subseteq I$, and let $\Lambda_2
\subseteq J$. Then the {\em cluster coherence} $\mu_c (\Lambda_1,
\Phi_1; \Phi_2)$ of $\Phi_1$ and $\Phi_2$ with respect to
$\Lambda_1$ is defined by
\[
\mu_c(\Lambda_1, \Phi_1; \Phi_2) = \max_{j \in J}  \sum_{i \in
\Lambda_1} | \langle \varphi_{1i}, \varphi_{2j} \rangle|,
\]
and the {\em cluster coherence} $\mu_c (\Phi_1; \Lambda_2, \Phi_2)$
of $\Phi_1$ and $\Phi_2$ with respect to $\Lambda_2$ is defined by
\[
\mu_c(\Phi_1; \Lambda_2, \Phi_2) = \max_{i \in I}  \sum_{j \in
\Lambda_2} | \langle \varphi_{1i}, \varphi_{2j} \rangle|.
\]
\end{definition}

The relation between joint concentration and cluster coherence is
made precise in the following result from \cite{DK10}.

\begin{proposition} [\cite{DK10}] \label{prop:kappamuc}
Let $\Phi_1 =(\varphi_{1i})_{i \in I}$ and $\Phi_2
=(\varphi_{2j})_{j \in J}$ be two Parseval frames for a Hilbert
space $\cH$, and let $\Lambda_1 \subseteq I$ and $\Lambda_2
\subseteq J$. Then
\[
\kappa(\Lambda_1, \Phi_1; \Lambda_2, \Phi_2) \le
\max\{\mu_c(\Lambda_1, \Phi_1; \Phi_2), \mu_c(\Phi_1; \Lambda_2,
\Phi_2)\}.
\]
\end{proposition}

\noindent {\bf Proof.}
Let $x \in \cH$. We now choose coefficient sequences $c_1$ and $c_2$
such that
\[
x=\Phi_1 c_1 = \Phi_2 c_2\] and, for $i=1,2$, \beq
\label{eq:smallestalpha} \norm{c_i}_1 \le \norm{d_i}_1 \quad
\mbox{for all }d_i \mbox{ with } x=\Phi_i d_i. \eeq This implies
that {\allowdisplaybreaks
\begin{eqnarray*}
\lefteqn{\norm{1_{\Lambda_1} \Phi_1^T x}_1 + \norm{1_{\Lambda_2} \Phi_2^T x}_1}\\[1ex]
& = & \norm{1_{\Lambda_1} \Phi_1^T \Phi_2 c_2}_1 + \norm{1_{\Lambda_2} \Phi_2^T \Phi_1 c_1}_1\\
& \le & \sum_{i \in \Lambda_1}\left(\sum_{j \in J}
\absip{\varphi_{1i}}{\varphi_{2j}} |c_{2j}|\right)
+ \sum_{j \in \Lambda_2}\left(\sum_{i \in I} \absip{\varphi_{1i}}{\varphi_{2j}} |c_{1i}|\right)\\
& = & \sum_{j \in J}\left(\sum_{i \in \Lambda_1}
\absip{\varphi_{1i}}{\varphi_{2j}}\right) |c_{2j}|
+ \sum_{i \in I} \left(\sum_{j \in \Lambda_2} \absip{\varphi_{1i}}{\varphi_{2j}}\right) |c_{1i}|\\[1ex]
& \le & \mu_c(\Lambda_1,\Phi_1;\Phi_2)
 \norm{c_2}_1 + \mu_c(\Lambda_2, \Phi_2; \Phi_1) \norm{c_1}_1\\[1ex]
& \le & \max\{ \mu_c(\Lambda_1,\Phi_1;\Phi_2),
\mu_c(\Lambda_2,\Phi_2;\Phi_1)\} (\norm{c_1}_1+\norm{c_2}_1).
\end{eqnarray*}
} Since $\Phi_1$ and $\Phi_2$ are Parseval frames, we have
\[
x = \Phi_i (\Phi_i^T \Phi_i c_i) \quad\mbox{for }i=1,2.
\]
Hence, by exploiting \eqref{eq:smallestalpha}, {\allowdisplaybreaks
\begin{eqnarray*}
\lefteqn{\norm{1_{\Lambda_1} \Phi_1^T x}_1 + \norm{1_{\Lambda_2} \Phi_2^T x}_1}\\[1ex]
& \le &\max\{ \mu_c(\Lambda_1,\Phi_1;\Phi_2), \mu_c(\Lambda_2,\Phi_2;\Phi_1)\} (\norm{\Phi_1^T \Phi_1 c_1}_1+\norm{\Phi_2^T \Phi_2 c_2}_1)\\[1ex]
& = & \max\{ \mu_c(\Lambda_1,\Phi_1;\Phi_2),
\mu_c(\Lambda_2,\Phi_2;\Phi_1)\} (\norm{\Phi_1^T x}_1+\norm{\Phi_2^T
x}_1). \quad  \mbox{ \qed}
\end{eqnarray*}
}

Before stating the data separation estimate which uses joint
concentration, we need to discuss the conditions on sparsity of the
components in the two Parseval frames. Since for real data `true
sparsity' is unrealistic, a weaker condition will be imposed.
For the next result, a notion invoking the clustering
of the significant coefficients will be required. This notion, first
utilized in \cite{Don06c}, is defined for our data separation
problem as follows.

\begin{definition}
Let $\Phi_1 =(\varphi_{1i})_{i \in I}$ and $\Phi_2
=(\varphi_{2j})_{j \in J}$ be two Parseval frames for a Hilbert
space $\cH$, and let $\Lambda_1 \subseteq I$ and $\Lambda_2
\subseteq J$. Further,  suppose that $x \in \cH$ can be decomposed
as $x=x_1^0+x_2^0$. Then the components $x_1^0$ and $x_2^0$ are
called {\em $\delta$-relatively sparse} in $\Phi_1$ and $\Phi_2$
with respect to $\Lambda_1$ and $\Lambda_2$, if
\[
\norm{1_{\Lambda_1^c} \Phi_1^T x_1^0}_1 + \norm{1_{\Lambda_2^c}
\Phi_2^T x_2^0}_1\le \delta.
\]
\end{definition}

We now have all ingredients to state the data separation result from
\cite{DK10}, which -- as compared to Theorem \ref{theo:generalbound}
-- now invokes information about the clustering of coefficients.

\begin{theorem} [\cite{DK10}] \label{theo:cluster}
Let $\Phi_1 =(\varphi_{1i})_{i \in I}$ and $\Phi_2
=(\varphi_{2j})_{j \in J}$ be two Parseval frames for a Hilbert
space $\cH$, and suppose that $x \in \cH$ can be decomposed as
$x=x_1^0+x_2^0$. Further, let $\Lambda_1 \subseteq I$ and $\Lambda_2
\subseteq J$ be chosen such that $x_1^0$ and $x_2^0$ are
$\delta$-relatively sparse in $\Phi_1$ and $\Phi_2$ with respect to
$\Lambda_1$ and $\Lambda_2$. Then the solution
$(x_1^\star,x_2^\star)$ of the $\ell_1$ minimization problem
$(\mbox{Sep}_{a})$ stated in \eqref{eq:newmingeneral} satisfies
\[
  \norm{x_1^\star-x_1^0}_2 + \norm{x_2^\star-x_2^0}_2 \le \frac{2\delta}{1-2\kappa}.
\]
\end{theorem}

\noindent {\bf Proof.}
First, using the fact that $\Phi_1$ and $\Phi_2$ are Parseval
frames,
\begin{eqnarray*}
\norm{x_1^\star-x_1^0}_2 + \norm{x_2^\star-x_2^0}_2
& = & \norm{\Phi_1^T(x_1^\star-x_1^0)}_2 + \norm{\Phi_2^T(x_2^\star-x_2^0)}_2\\
& \le & \norm{\Phi_1^T(x_1^\star-x_1^0)}_1 +
\norm{\Phi_2^T(x_2^\star-x_2^0)}_1.
\end{eqnarray*}
The decomposition $x_1^0+x_2^0=x=x_1^\star+x_2^\star$ implies
\[
x_2^\star-x_2^0 = -(x_1^\star-x_1^0),
\]
which allows us to conclude that \beq \label{eq:gensepres1}
\norm{x_1^\star-x_1^0}_2 + \norm{x_2^\star-x_2^0}_2 \le
\norm{\Phi_1^T(x_1^\star-x_1^0)}_1 +
\norm{\Phi_2^T(x_1^\star-x_1^0)}_1. \eeq By the definition of
$\kappa$,
\begin{eqnarray*}
\lefteqn{\norm{\Phi_1^T(x_1^\star-x_1^0)}_1 + \norm{\Phi_2^T(x_1^\star-x_1^0)}_1}\\
& = & (\norm{1_{\Lambda_1}\Phi_1^T(x_1^\star-x_1^0)}_1 +
\norm{1_{\Lambda_2}\Phi_2^T(x_1^\star-x_1^0)}_1)  +
\norm{1_{\Lambda_1^c}\Phi_1^T(x_1^\star-x_1^0)}_1\\
& &  + \norm{1_{\Lambda_2^c}\Phi_2^T(x_2^\star-x_2^0)}_1 \\
& \le & \kappa \cdot \left( \norm{\Phi_1^T(x_1^\star-x_1^0)}_1 +
\norm{\Phi_2^T(x_1^\star-x_1^0)}_1 \right) +
\norm{1_{\Lambda_1^c}\Phi_1^T(x_1^\star-x_1^0)}_1\\
& &  + \norm{1_{\Lambda_2^c}\Phi_2^T(x_2^\star-x_2^0)}_1,
\end{eqnarray*}
which yields
\begin{eqnarray} \nonumber
\lefteqn{\norm{\Phi_1^T(x_1^\star-x_1^0)}_1 +
\norm{\Phi_2^T(x_1^\star-x_1^0)}_1}\\ \nonumber & \le &
\frac{1}{1-\kappa}
(\norm{1_{\Lambda_1^c}\Phi_1^T(x_1^\star-x_1^0)}_1 +
\norm{1_{\Lambda_2^c}\Phi_2^T(x_2^\star-x_2^0)}_1)\\ \nonumber & \le
& \frac{1}{1-\kappa} (\norm{1_{\Lambda_1^c}\Phi_1^Tx_1^\star}_1 +
\norm{1_{\Lambda_1^c}\Phi_1^Tx_1^0}_1 +
\norm{1_{\Lambda_2^c}\Phi_2^Tx_2^\star}_1+\norm{1_{\Lambda_2^c}\Phi_2^Tx_2^0}_1).
\end{eqnarray}
Now using the relative sparsity of $x_1^0$ and $x_2^0$ in $\Phi_1$
and $\Phi_2$ with respect to $\Lambda_1$ and $\Lambda_2$, we obtain
\beq \label{eq:gensepres2} \norm{\Phi_1^T(x_1^\star-x_1^0)}_1 +
\norm{\Phi_2^T(x_1^\star-x_1^0)}_1 \le \frac{1}{1-\kappa}
\left(\norm{1_{\Lambda_1^c}\Phi_1^Tx_1^\star}_1+
\norm{1_{\Lambda_2^c}\Phi_2^Tx_2^\star}_1+ \delta \right). \eeq By
the minimality of $x_1^\star$ and $x_2^\star$ as solutions of
$(\mbox{Sep}_{a})$ implying that
\begin{eqnarray*}
\sum_{i=1}^2
\left(\norm{1_{\Lambda_i^c}\Phi_i^Tx_i^\star}_1+\norm{1_{\Lambda_i}\Phi_i^Tx_i^\star}_1\right)
& = & \norm{\Phi_1^Tx_1^\star}_1+\norm{\Phi_2^Tx_2^\star}_1\\
& \le & \norm{\Phi_1^T x_1^0}_1+\norm{\Phi_2^T x_2^0}_1,
\end{eqnarray*}
we have
\begin{eqnarray*}
\lefteqn{\norm{1_{\Lambda_1^c}\Phi_1^Tx_1^\star}_1+ \norm{1_{\Lambda_2^c}\Phi_2^Tx_2^\star}_1}\\
& \le & \norm{\Phi_1^T x_1^0}_1+\norm{\Phi_2^T x_2^0}_1 - \norm{1_{\Lambda_1}\Phi_1^Tx_1^\star}_1 - \norm{1_{\Lambda_2}\Phi_2^Tx_2^\star}_1\\
& \le & \norm{\Phi_1^T x_1^0}_1+\norm{\Phi_2^T x_2^0}_1
+ \norm{1_{\Lambda_1}\Phi_1^T(x_1^\star-x_1^0)}_1 - \norm{1_{\Lambda_1}\Phi_1^T x_1^0}_1\\
& & + \norm{1_{\Lambda_2}\Phi_2^T(x_2^\star-x_2^0)}_1 -
\norm{1_{\Lambda_2}\Phi_2^T x_2^0}_1.
\end{eqnarray*}
Again exploiting relative sparsity leads to \beq
\label{eq:gensepres3} \norm{1_{\Lambda_1^c}\Phi_1^Tx_1^\star}_1+
\norm{1_{\Lambda_2^c}\Phi_2^Tx_2^\star}_1 \le
\norm{1_{\Lambda_1}\Phi_1^T(x_1^\star-x_1^0)}_1 +
\norm{1_{\Lambda_2}\Phi_2^T(x_2^\star-x_2^0)}_1 + \delta. \eeq
Combining \eqref{eq:gensepres2} and \eqref{eq:gensepres3} and again
using joint concentration,
\begin{eqnarray*}
\lefteqn{\norm{\Phi_1^T(x_1^\star-x_1^0)}_1 + \norm{\Phi_2^T(x_1^\star-x_1^0)}_1}\\
& \le & \frac{1}{1-\kappa} \left
[\norm{1_{\Lambda_1}\Phi_1^T(x_1^\star-x_1^0)}_1 +
\norm{1_{\Lambda_2}\Phi_2^T(x_1^\star-x_1^0)}_1+
2\delta \right]  \\
& \le & \frac{1}{1-\kappa} \left[\kappa \cdot
(\norm{\Phi_1^T(x_1^\star-x_1^0)}_1 +
\norm{\Phi_2^T(x_1^\star-x_1^0)}_1)+ 2\delta \right].
\end{eqnarray*}
Thus, by \eqref{eq:gensepres1}, we finally obtain
\[
\norm{x_1^\star-x_1^0}_2 + \norm{x_2^\star-x_2^0}_2 \le \left(1-
\frac{\kappa}{1-\kappa}\right)^{-1} \cdot \frac{2\delta}{1-\kappa} =
\frac{2\delta}{1-2\kappa}. \quad \mbox{ \qed}
\]

Using Proposition \ref{prop:kappamuc}, this result can also be
stated in terms of cluster coherence, which on the one hand provides
an easier accessible estimate and allows a better comparison with
results using mutual coherence, but on the other hand poses a
slightly weaker estimate.

\begin{theorem} [\cite{DK10}] \label{theo:cluster2}
Let $\Phi_1 =(\varphi_{1i})_{i \in I}$ and $\Phi_2
=(\varphi_{2j})_{j \in J}$ be two Parseval frames for a Hilbert
space $\cH$, and suppose that $x \in \cH$ can be decomposed as
$x=x_1^0+x_2^0$. Further, let $\Lambda_1 \subseteq I$ and $\Lambda_2
\subseteq J$ be chosen such that $x_1^0$ and $x_2^0$ are
$\delta$-relatively sparse in $\Phi_1$ and $\Phi_2$ with respect to
$\Lambda_1$ and $\Lambda_2$. Then the solution
$(x_1^\star,x_2^\star)$ of the minimization problem
$(\mbox{Sep}_{a})$ stated in \eqref{eq:newmingeneral} satisfies
\[
  \norm{x_1^\star-x_1^0}_2 + \norm{x_2^\star-x_2^0}_2 \le \frac{2\delta}{1-2\mu_c},
\]
with
\[
\mu_c = \max\{\mu_c(\Lambda_1, \Phi_1; \Phi_2), \mu_c(\Phi_1;
\Lambda_2, \Phi_2)\}.
\]
\end{theorem}

To thoroughly understand this estimate, it is important to notice
that both relative sparsity $\delta$ as well as cluster coherence
$\mu_c$ depend heavily on the choice of the sets of significant
coefficients $\Lambda_1$ and $\Lambda_2$. Choosing those sets too
large allows for a very small $\delta$, however $\mu_c$ might not be
less than $\frac12$ anymore, thereby making the estimate useless.
Choosing those sets too small will force $\mu_c$ to become
simultaneously small, in particular, smaller than $\frac12$, with
the downside that $\delta$ might be large.

It is also essential to realize that the sets $\Lambda_1$ and
$\Lambda_2$ are a mere analysis tool; they do not appear in the
minimization problem $(\mbox{Sep}_{a})$. This means that the
algorithm does not care about this choice at all, however the
estimate for accuracy of separation does.

Also note that this result can be easily generalized to general
frames instead of Parseval frames, which then changes the separation
estimate by invoking the lower frame bound. In addition, a version
including noise was derived in \cite{DK10}.


\subsection{Relation with Uncertainty Principles}

Intriguingly, there exists a very close connection between
uncertainty principles and data separation problems. Given a signal
$x \in \cH$ and two bases or frames $\Phi_1$ and $\Phi_2$, loosely
speaking, an uncertainty principle states that $x$ cannot be
sparsely represented by $\Phi_1$ and $\Phi_2$ simultaneously; one of
the expansions is always not sparse unless $x=0$. For the relation
to the `classical' uncertainty principle, we refer to Subsection
\ref{subsec:sinspike}.

The first result making this uncertainty viewpoint precise was
proven in \cite{EB02} with ideas already lurking in \cite{DS89} and
\cite{DH01}. Again, it turns out that the mutual coherence is an
appropriate measure for allowed sparsity, here serving as a lower
bound for the simultaneously achievable sparsity of two expansions.

\begin{theorem} [\cite{EB02}] \label{theo:UCP1}
Let $\Phi_1$ and $\Phi_2$ be two orthonormal bases for a Hilbert
space $\cH$, and let $x \in \cH$, $x \neq 0$. Then
\[
\|\Phi_1^T x\|_0 + \|\Phi_2^T x\|_0 \ge
\frac{2}{\mu([\Phi_1|\Phi_2])}.
\]
\end{theorem}

\noindent {\bf Proof.}
First, let $\Phi_1 =(\varphi_{1i})_{i \in I}$ and $\Phi_2
=(\varphi_{2j})_{j \in J}$. Further, let $\Lambda_1 \subseteq I$ and
$\Lambda_2 \subset J$ denote the support of $\Phi_1^T x$ and
$\Phi_2^T x$, respectively. Since $x = \Phi_1 \Phi_1^T x$, for each
$j \in J$, \beq \label{eq:UCP11} |(\Phi_2^T x)_j| = \left|\sum_{i
\in \Lambda_1} (\Phi_1^T x)_i
\ip{\varphi_{1i}}{\varphi_{2j}}\right|. \eeq Since $\Phi_1$ and
$\Phi_2$ are orthonormal bases, we have \beq \label{eq:UCP12}
\|x\|_2 = \|\Phi_1^T x\|_2 = \|\Phi_2^T x\|_2. \eeq Using in
addition the Cauchy-Schwarz inequality, we can continue
\eqref{eq:UCP11} by
\[
|(\Phi_2^T x)_j|^2 \le \norm{\Phi_1^T x}_2^2 \cdot \left|\sum_{i \in
\Lambda_1} \absip{\varphi_{1i}}{\varphi_{2j}}^2 \right| \le
\|x\|_2^2 \cdot |\Lambda_1| \cdot \mu([\Phi_1|\Phi_2])^2.
\]
This implies
\[
\|\Phi_2^T x\|_2 = \left(\sum_{j \in \Lambda_2} |(\Phi_2^T
x)_j|^2\right)^{1/2}  \le \|x\|_2 \cdot \sqrt{|\Lambda_1| \cdot
|\Lambda_2|} \cdot \mu([\Phi_1|\Phi_2]).
\]
Since $|\Lambda_i| = \|\Phi_i^T x\|_0$, $i = 1, 2$, and again using
\eqref{eq:UCP12}, we obtain
\[
\sqrt{\|\Phi_1^T x\|_0  \cdot \|\Phi_2^T x\|_0} \ge
\frac{1}{\mu([\Phi_1|\Phi_2])}.
\]
Using the geometric-algebraic relationship,
\[
\frac{1}{2} (\|\Phi_1^T x\|_0 + \|\Phi_2^T x\|_0) \ge
\sqrt{\|\Phi_1^T x\|_0  \cdot \|\Phi_2^T x\|_0} \ge
\frac{1}{\mu([\Phi_1|\Phi_2])},
\]
which proves the claim.
\qed

This result can be easily connected to the problem of simultaneously
sparse expansions. The following version was first explicitly stated
in \cite{BDE09}.

\begin{theorem} [\cite{BDE09}] \label{theo:UCP2}
Let $\Phi_1$ and $\Phi_2$ be two orthonormal bases for a Hilbert
space $\cH$, and let $x \in \cH$, $x \neq 0$. Then, for any two
distinct coefficient sequences $c_i$ satisfying $x =
[\Phi_1|\Phi_2]c_i$, $i=1, 2$, we have
\[
\|c_1\|_0 + \|c_2\|_0 \ge \frac{2}{\mu([\Phi_1|\Phi_2])},
\]
\end{theorem}

\noindent {\bf Proof.}
First, set $d = c_1-c_2$ and partition $d$ into
$[d_{\Phi_1},d_{\Phi_2}]^T$ such that
\[
0 = [\Phi_1|\Phi_2]d = \Phi_1 d_{\Phi_1} + \Phi_2 d_{\Phi_2}.
\]
Since $\Phi_1$ and $\Phi_2$ are bases and $d \neq 0$, the vector $y$
defined by
\[
y = \Phi_1 d_{\Phi_1} = - \Phi_2 d_{\Phi_2}
\]
is non-zero. Applying Theorem \ref{theo:UCP1}, we obtain
\[
\|d\|_0 = \|d_{\Phi_1}\|_0 + \|d_{\Phi_2}\|_0 \ge
\frac{2}{\mu([\Phi_1|\Phi_2])}.
\]
Since $d = c_1-c_2$, we have
\[
\|c_1\|_0 + \|c_2\|_0 \ge \|d\|_0 \ge
\frac{2}{\mu([\Phi_1|\Phi_2])}. \quad \mbox{ \qed}
\]

We would also like to mention the very recent paper \cite{Tro10} by
Tropp, in which he studies uncertainty principles for random sparse
signals over an incoherent dictionary. He, in particular, shows that
the coefficient sequence of each non-optimal expansion of a signal
contains far more non-zero entries than the one of the sparsest
expansion.


\section{Signal Separation}
\label{sec:signals}

In this section, we study the special situation of signal
separation, where we refer to 1D signals as opposed to images, etc.
For this, we start with the most prominent example of separating
sinusoids from spikes, and then discuss further problem classes.


\subsection{Separation of Sinusoids and Spikes}
\label{subsec:sinspike}

Sinusoidal and spike components are intuitively the morphologically
most distinct features of a signal, since one is periodic and the
other transient. Thus, it seems natural that the first results using
sparsity and $\ell_1$ minimization for data separation were proven
for this situation. Certainly, real-world signals are never a
pristine combination of sinusoids and spikes. However, thinking of
audio data from a recording of musical instruments, these components
are indeed an essential part of such signals.

The separation problem can be generally stated in the following way:
Let the vector $x \in \mathbb{R}^n$ consist of $n$ samples of a
continuum domain signal at times $t \in \{0, \ldots, n-1\}$. We
assume that $x$ can be decomposed into
\[
x = x_1 + x_2.
\]

Here $x_1$ shall consist of  $n$ samples -- at the same points in
time as $x$ -- of a continuum domain signal of the form
\[
\frac{1}{\sqrt{n}} \sum_{\omega = 0}^{n-1} c_{1\omega} e^{2 \pi i
\omega t/n}, \quad t \in \mathbb{R}.
\]
Thus, by letting $\Phi_1 = (\varphi_{1\omega})_{0 \le \omega \le
n-1}$ denote the Fourier basis, i.e.,
\[
\varphi_{1\omega} = \left(\tfrac{1}{\sqrt{n}} e^{2 \pi i \omega
t/n}\right)_{0 \le t \le n-1},
\]
the discrete signal $x_1$ can be written as
\[
x_1 = \Phi_1 c_1 \quad \mbox{with } c_1 = (c_{1\omega})_{0 \le
\omega \le n-1}.
\]
If $x_1$ is now the superposition of very few sinusoids, then the
coefficient vector $c_1$ is sparse.

Further, consider a continuum domain signal which has a few spikes.
Sampling this signal at $n$ samples at times $t \in \{0, \ldots,
n-1\}$ leads to a discrete signal $x_2 \in \mathbb{R}^n$ which has
very few non-zero entries. In order to expand $x_2$ in terms of a
suitable representation system, we let $\Phi_2$ denote the Dirac
basis, i.e., $\Phi_2$ is simply the identity matrix, and write
\[
x_2 = \Phi_2 c_2,
\]
where $c_2$ is then a sparse coefficient vector.

\begin{figure}[h]
\centering

\includegraphics[scale=0.21]{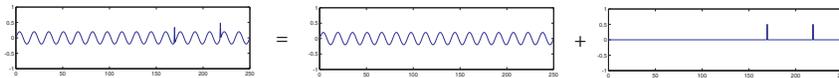}

\vspace*{0.25cm}

\caption{Separation of artificial audio data into sinusoids and
spikes.} \label{fig:sinusoidspike}
\end{figure}

The task now consists in extracting $x_1$ and $x_2$ from the known
signal $x$, which is illustrated in Figure \ref{fig:sinusoidspike}.
It will be illuminating to detect the dependence on the number of
sampling points of the bound for the sparsity of $c_1$ and $c_2$
which still allows for separation via $\ell_1$ minimization.

The intuition that -- from a morphological standpoint -- this
situation is extreme, can be seen by computing the mutual coherence
between the Fourier basis $\Phi_1$ and the Dirac basis $\Phi_2$. For
this, we obtain \beq \label{eq:muFD} \mu([\Phi_1|\Phi_2]) =
\frac{1}{\sqrt{n}}, \eeq and, in fact, $1/\sqrt{n}$ is the minimal
possible value. This can be easily seen: If $\Phi_1$ and $\Phi_2$
are two general orthonormal bases of $\mathbb{R}^n$, then $\Phi_1^T
\Phi_2$ is an orthonormal matrix. Hence the sum of squares of its
entries equals $n$, which implies that all entries can not be less
than $1/\sqrt{n}$.

The following result from \cite{EB02} makes this dependence precise.
We wish to mention that the first answer to this question was
derived in \cite{DH01}. In this paper the slightly weaker bound of
$(1+\sqrt{n})/2$ for $\|c_1\|_0 + \|c_2\|_0$ was proven by using the
general result in Theorem \ref{theo:generalbound} instead of the
more specialized Theorem \ref{theo:specialbound} exploited to derive
the result from \cite{EB02} stated below.

\begin{theorem} [\cite{EB02}]
Let $\Phi_1$ be the Fourier basis for $\mathbb{R}^n$ and let
$\Phi_2$ be the Dirac basis for $\mathbb{R}^n$. Further, let $x \in
\mathbb{R}^n$ be the signal
\[
x = x_1 + x_2, \quad \mbox{where } x_1 = \Phi_1 c_1 \mbox{ and } x_2
= \Phi_2 c_2,
\]
with coefficient vectors $c_i \in \mathbb{R}^n$, $i= 1, 2$. If
\[
\|c_1\|_0 + \|c_2\|_0 < (\sqrt{2}-0.5)\sqrt{n},
\]
then the $\ell_1$ minimization problem $(\mbox{Sep}_{s})$ stated in
\eqref{eq:mingeneral} recovers $c_1$ and $c_2$ uniquely, and hence
extracts $x_1$ and $x_2$ from $x$ precisely.
\end{theorem}

\noindent {\bf Proof.}
Recall that we have (cf. \eqref{eq:muFD})
\[
\mu([\Phi_1|\Phi_2]) = \frac{1}{\sqrt{n}}.
\]
Hence, by Theorem \ref{theo:specialbound}, the $\ell_1$ minimization
problem $(\mbox{Sep}_{s})$ recovers $c_1$ and $c_2$ uniquely,
provided that
\[
\|c_1\|_0 + \|c_2\|_0 < \frac{\sqrt{2}-0.5}{\mu([\Phi_1|\Phi_2])} =
(\sqrt{2}-0.5)\sqrt{n}.
\]
The theorem is proved.
\qed

The classical uncertainty principle states that, roughly speaking, a
function cannot both be localized in time as well as in frequency
domain. A discrete version of this fundamental principle was --
besides the by now well-known continuum domain Donoho-Stark
uncertainty principle -- derived in \cite{DS89}. It showed that a
discrete signal and its Fourier transform cannot both be highly
localized in the sense of having `very few' non-zero entries. We
will now show that this result -- as it was done in \cite{DH01} --
can be interpreted as a corollary from data separation results.

\begin{theorem} [\cite{DS89}]
Let $x \in \mathbb{R}^n$, and denote its Fourier transform by
$\hat{x}$. Then
\[
\|x\|_0 + \|\hat{x}\|_0 \ge 2 \sqrt{n}.
\]
\end{theorem}

\noindent {\bf Proof.}
For the proof, we intend to use Theorem \ref{theo:UCP1}. First, we
note that by letting $\Phi_1$ denote the Dirac basis, we trivially
have
\[
\|\Phi_1^T x\|_0 = \|x\|_0.
\]
Secondly, letting $\Phi_2$ denote the Fourier basis, we obtain
\[
\hat{x} = \Phi_2^T x.
\]
Now recalling that, by \eqref{eq:muFD},
\[
\mu([\Phi_1|\Phi_2]) = \frac{1}{\sqrt{n}},
\]
we can conclude from Theorem \ref{theo:UCP1} that
\[
\|x\|_0 + \|\hat{x}\|_0 = \|\Phi_1^T x\|_0 + \|\Phi_2^T x\|_0 \ge
\frac{2}{\mu([\Phi_1|\Phi_2])} = 2 \sqrt{n}.
\]
This finishes the proof.
\qed

As an excellent survey about sparsity of expansions of signals in
the Fourier and Dirac basis, data separation, and related
uncertainty principles as well as on very recent results using
random signals, we refer to \cite{Tro08}.


\subsection{Further Variations}

Let us briefly mention the variety of modifications of the previous
discussed setting, most of them empirical analyses, which were
developed during the last few years.

The most common variation of the sinusoid and spike setting is the
consideration of a more general periodic component, which is then
considered to be sparse in a Gabor system, superimposed by a second
component, which is considered to be sparse in a system sensitive to
spike-like structures similar to wavelets. This is, for instance,
the situation considered in \cite{GB03}. An example for a different
setting is the substitution of a Gabor system by a Wilson basis,
analyzed in \cite{BGN08}. In this paper, as already mentioned in
Subsection \ref{subsec:clusteredsparsity}, the clustering of
coefficients already plays an essential role. It should also be
mentioned that a specifically adapted norm, namely the mixed
$\ell_{1,2}$ or $\ell_{2,1}$ norm, is used in \cite{KT09} to take
advantage of this clustering, and various numerical experiments show
successful separation.


\section{Image Separation}
\label{sec:images}

This section is devoted to discuss results on image separation
exploiting Morphological Component Analysis, first focussing on
empirical studies and secondly on theoretical results.


\subsection{Empirical Results}
\label{subsec:empiricalimages}

In practice, the observed signal $x$ is often contaminated by noise,
i.e., $x = x_1 + x_2 + n$ containing the to-be-extracted components
$x_1$ and $x_2$ and some noise $n$. This requires an adaption of the
$\ell_1$ minimization problem. As proposed in numerous publications,
one typically considers a modified optimization problem -- so-called
{\em Basis Pursuit Denoising} -- which can be obtained by relaxing
the constraint in order to deal with noisy observed signals. The
$\ell_1$ minimization problem $(\mbox{Sep}_{s})$ stated in
\eqref{eq:mingeneral}, which places the $\ell_1$ norm on the {\em
synthesis} side then takes the form:
\[
\min_{c_1,c_2} \norm{c_1}_1 + \norm{c_2}_1 + \lambda\|x - \Phi_1 c_1
- \Phi_2 c_2\|_2^2
\]
with appropriately chosen regularization parameter $\lambda > 0$.
Similarly, we can consider the relaxed form of the $\ell_1$
minimization problem $(\mbox{Sep}_{a})$ stated in
\eqref{eq:newmingeneral}, which places the $\ell_1$ norm on the {\em
analysis} side:
\[
\min_{x_1, x_2} \norm{\Phi_1^T x_1}_1 + \norm{\Phi_2^T x_2}_1 +
\lambda\|x - x_1 - x_2\|_2^2.
\]
In these new forms, the additional content in the image -- the noise
--, characterized by the property that it can not be represented
sparsely by either one of the two systems $\Phi_1$ and $\Phi_2$,
will be allocated to the residual $x - \Phi_1 c_1 - \Phi_2 c_2$ or
$x - x_1 - x_2$ depending on which of the two minimization problems
stated above is chosen. Hence, performing this minimization, we not
only separate the data, but also succeed in removing an additive
noise component as a by-product.

There exist by now a variety of algorithms which numerically solve
such minimization problems. One large class are, for instance,
iterative shrinkage algorithms; and we refer to the beautiful new
book \cite{Ela10} by Elad for an overview. It should be mentioned
that it is also possible to perform these separation procedures
locally, thus enabling parallel processing, and again we refer to
\cite{Ela10} for further details.

\medskip

Let us now delve into more concrete situations. One prominent class
of empirical studies concerns the separation of point- and curvelike
structures. This type of problem arises, for instance, in
astronomical imaging, where astronomers would like to separate stars
(pointlike structures) from filaments (curvelike structures).
Another area in which the separation of points from curves is
essential is neurobiological imaging. In particular, for Alzheimer
research, neurobiologists analyze images of neurons, which --
considered in 2D -- are a composition of the dendrites (curvelike
structures) of the neuron and the attached spines (pointlike
structures). For further analysis of the shape of these components,
dendrites and spines need to be separated.

From a mathematical perspective, pointlike structures are generally
speaking 0D structures whereas curvelike structures are 1D
structures, which reveals their morphological difference. Thus it
seems conceivable that separation using the idea of Morphological
Component Analysis can be achieved, and the empirical results
presented in the sequel as well as the theoretical results discussed
in Subsection \ref{subsec:theoryimages} give evidence to this claim.

To set up the minimization problem properly, the question arises
which systems adapted to the point- and curvelike objects to use.
For extracting pointlike structures, wavelets seem to be optimal,
since they provide optimally sparse approximations of smooth
functions with finitely many point singularities. As a sparsifying
system for curvelike structures, two different possibilities were
explored so far. From a historical perspective, the first system to
be utilized were {\em curvelets} \cite{CD05}, which provide
optimally sparse approximations of smooth functions exhibiting
curvilinear singularities. The composed dictionary of
wavelets-curvelets is used in MCALab\footnote{MCALab (Version 120)
is available from \url{http://jstarck.free.fr/jstarck/Home.html}.},
and implementation details are provided in the by now considered
fundamental paper \cite{SED05}.
A few years later {\em shearlets} were developed, see \cite{GKL06}
or the survey paper \cite{KLL10}, which deal with curvilinear
singularities in a similarly favorable way as curvelets (cf.
\cite{KL10a}), but have, for instance, the advantage of providing a
unified treatment of the continuum and digital realm and being
associated with a fast transform. Separation using the resulting
dictionary of wavelets-shearlets is implemented and publicly
available in ShearLab\footnote{ShearLab (Version 1.1) is available
from \url{http://www.shearlab.org}.}. For a close comparison between
both approaches we refer to \cite{KL10b} -- in this paper the
separation algorithm using wavelets and shearlets is also detailed
--, where a numerical comparison shows that ShearLab provides a
faster as well as more precise separation.

\begin{figure}[ht!]
\centering
\subfloat[Original image]{%
\includegraphics[width=4.7cm, height=4.5cm]{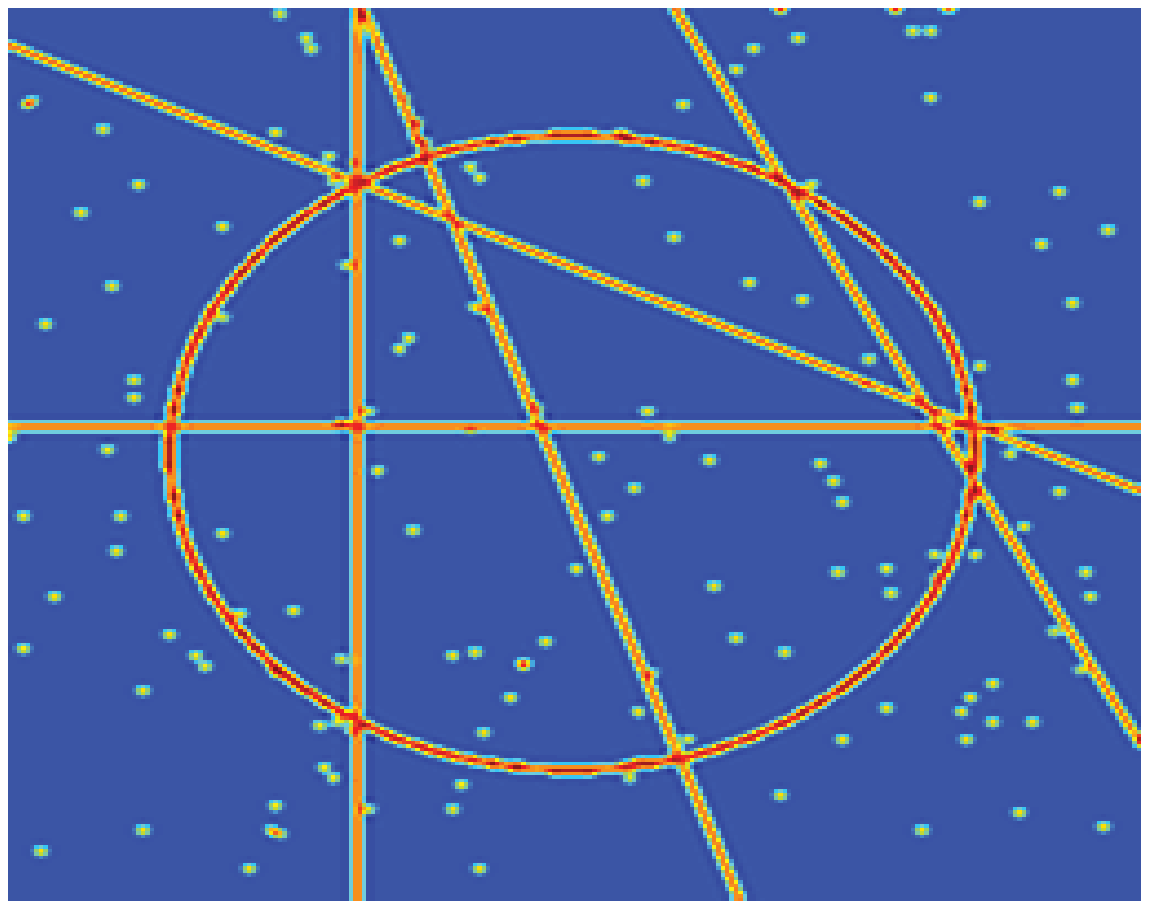}
\label{fig:curvepoint1}}
\subfloat[Noisy image]{%
\includegraphics[width=4.7cm, height=4.5cm]{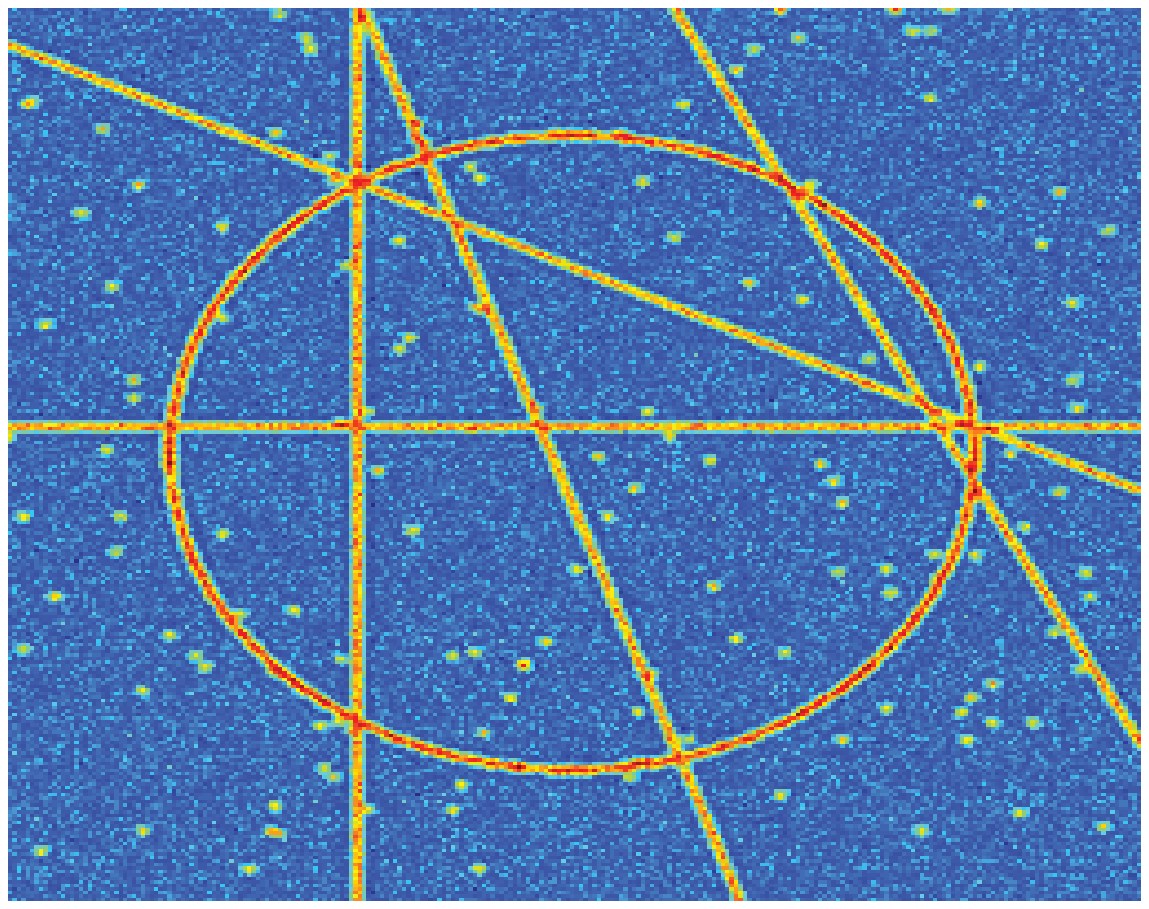}
\label{fig:noisy}}

\subfloat[Pointlike Component]{%
\includegraphics[width=4.7cm, height=4.5cm]{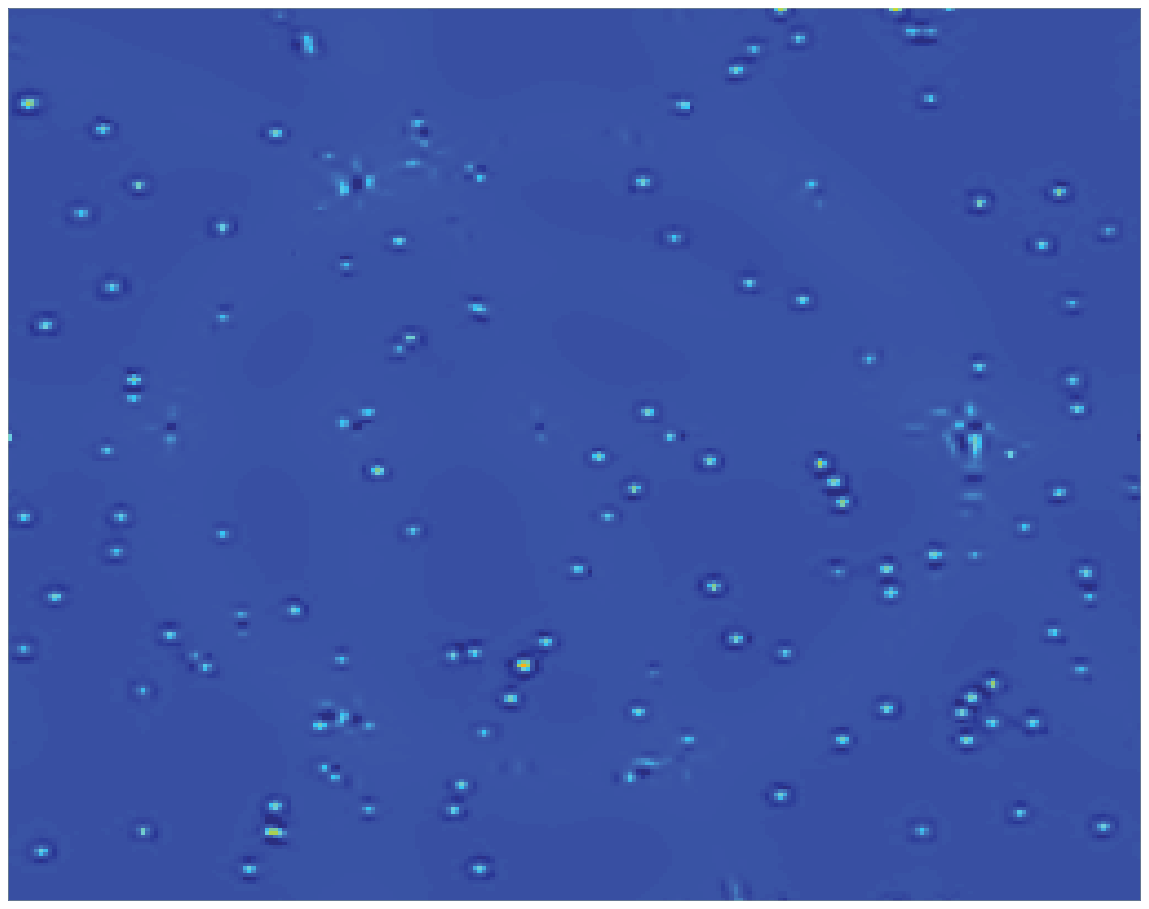}
\label{fig:shearpoint1} }
\subfloat[Curvelike Component]{%
\includegraphics[width=4.7cm, height=4.5cm]{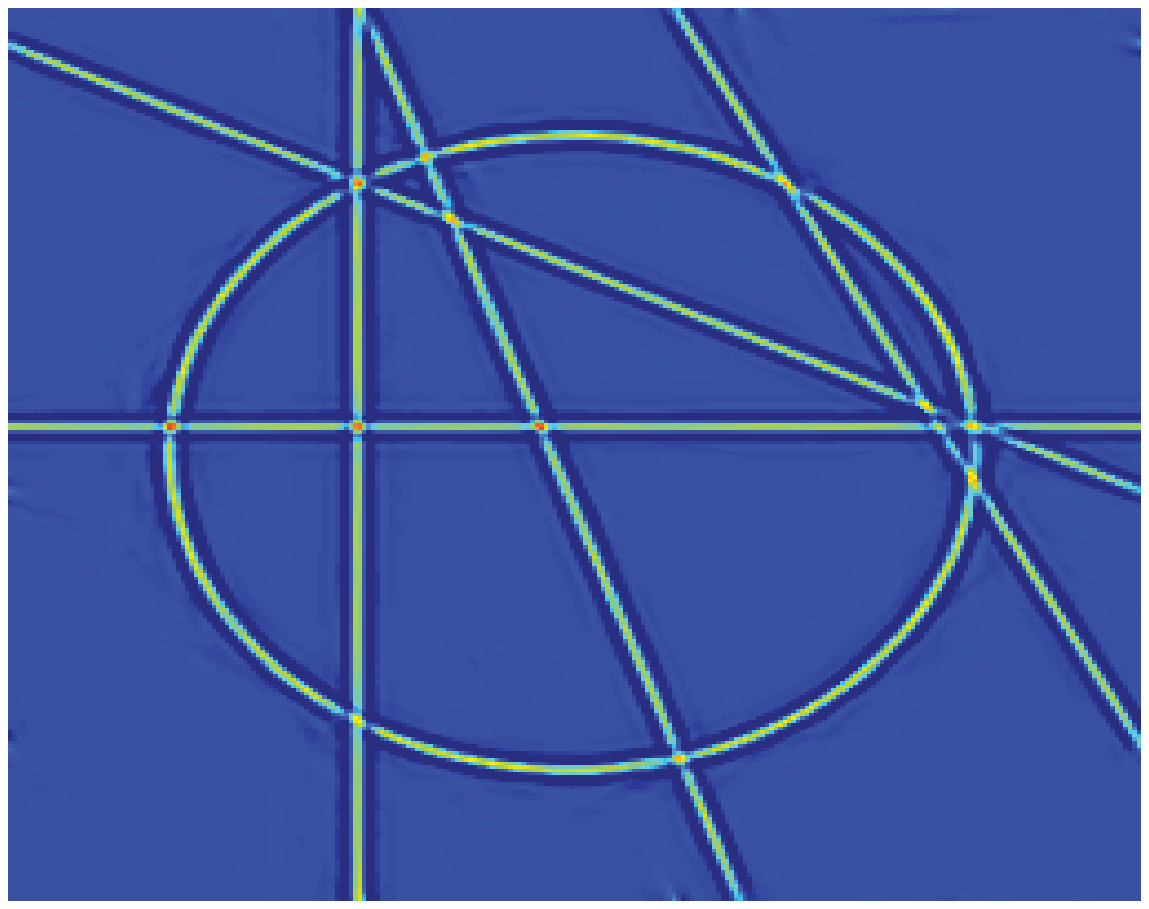}
\label{fig:shearcurve1} }

\vspace*{0.25cm}

\caption{Separation of an artificial image composed of points,
lines, and a circle into point- and curvelike components using
ShearLab.} \label{fig:circle}
\end{figure}

\begin{figure}[ht!]
\centering

\subfloat[Original image]{%
\includegraphics[scale=0.4]{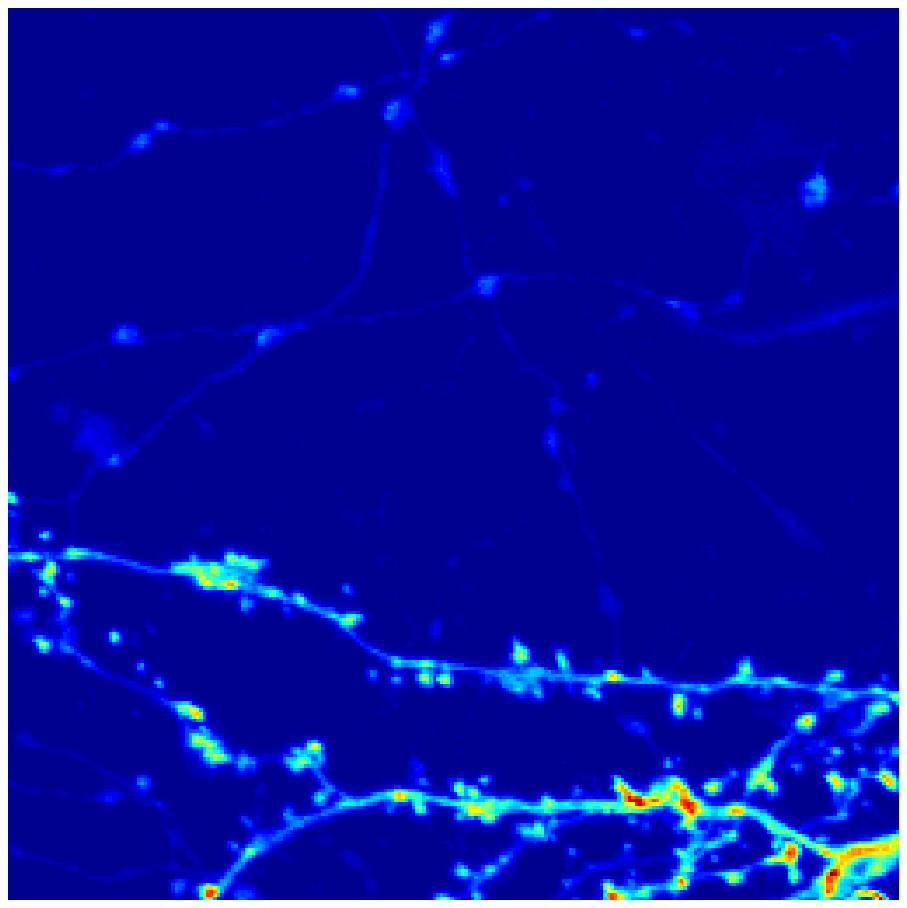}
\label{fig:curvepoint2}}

\subfloat[Pointlike Component]{%
\includegraphics[scale=0.4]{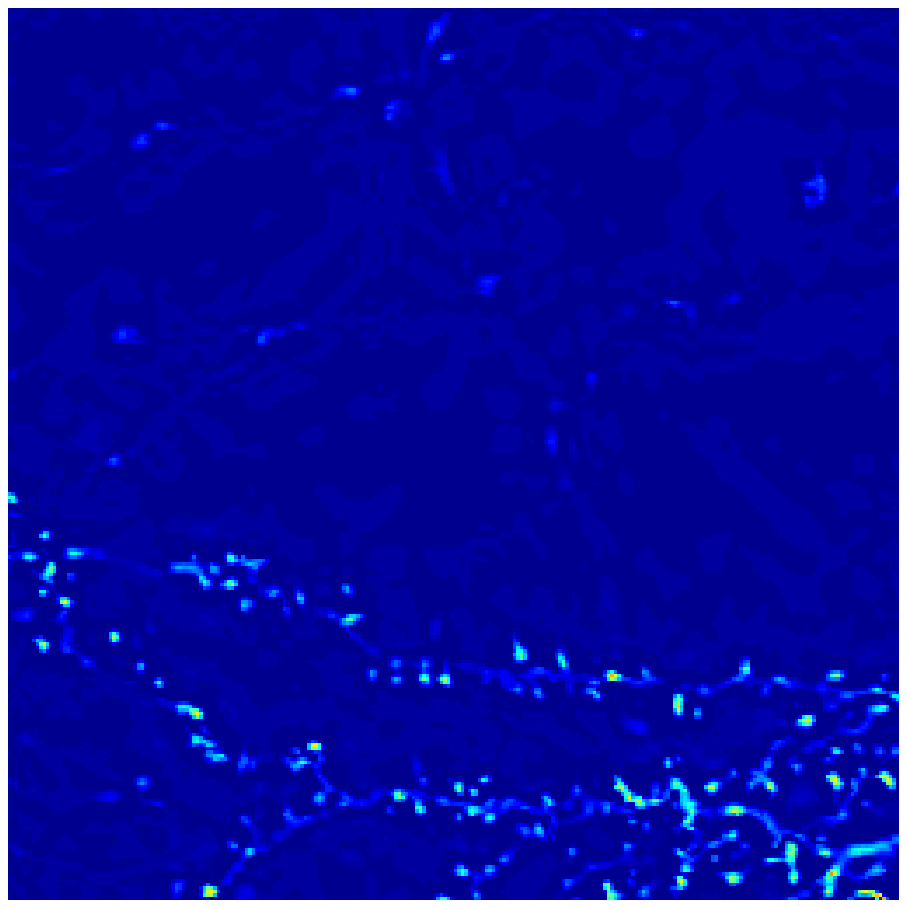}
\label{fig:shearpoint2} }
\subfloat[Curvelike Component]{%
\includegraphics[scale=0.4]{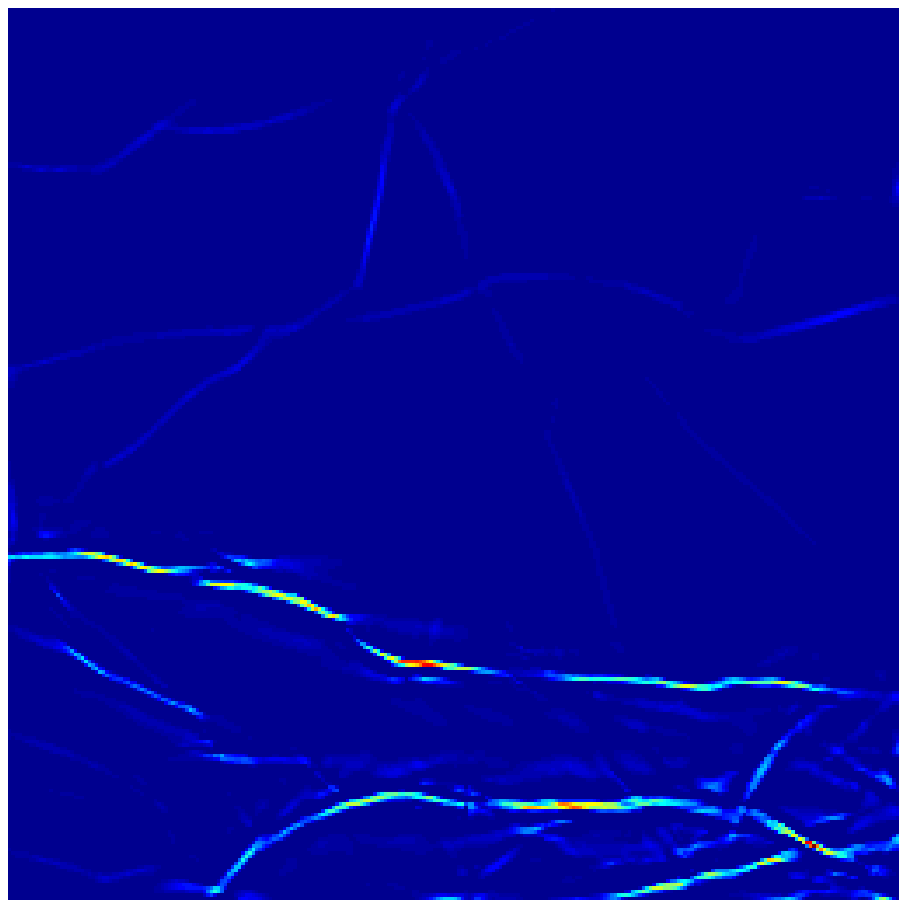}
\label{fig:shearcurve2} }

\vspace*{0.25cm}

\caption{Separation of a neuron image into point- and curvelike
components using ShearLab.} \label{fig:neuron}
\end{figure}

For illustrative purposes, Figure \ref{fig:circle} shows the
separation of an artificial image composed of points, lines, and a
circle as well as added noise into the pointlike structures (points)
and the curvelike structures (lines and the circle), while removing
the noise simultaneously. The only visible artifacts can be seen at
the intersections of the curvelike structures, which is not
surprising since it is even justifiable to label these intersections
as `points'. As an example using real data, we present in Figure
\ref{fig:neuron} the separation of a neuron image into dendrites and
spines again using ShearLab.

\medskip

\begin{figure}[ht!]
\centering

\includegraphics[height=1in]{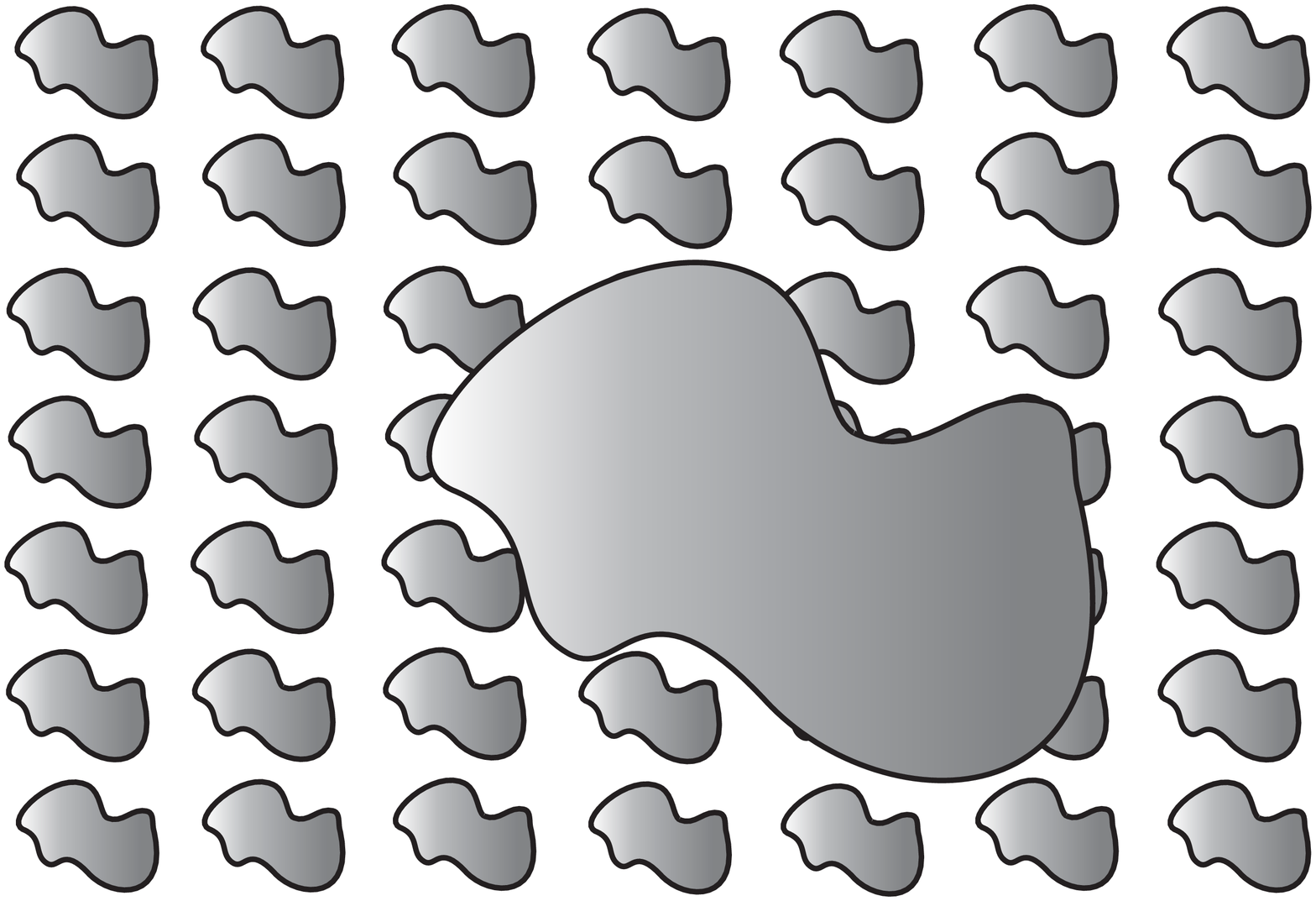}

\vspace*{0.25cm}

\caption{Periodic small cartoons versus one large cartoon.}
\label{fig:periodiccartoons}
\end{figure}

Another widely explored category of image separation is the
separation of cartoons and texture. Here, the term cartoon typically
refers to a piecewise smooth part in the image, and texture means a
periodic structure. A mathematical model for a {\em cartoon} was
first introduced in \cite{Don01} as a $C^2$ function containing a
$C^2$ discontinuity. In contrast to this, the term {\em texture} is
a widely open expression, and people have debated for years over an
appropriate model for the texture content of an image. A viewpoint
from applied harmonic analysis characterizes texture as a structure
which exhibits a sparse expansion in a Gabor system. As a side
remark, the reader should be aware that periodizing a cartoon part
of an image produces a texture component, thereby revealing the very
fine line between cartoons and texture, illustrated in Figure
\ref{fig:periodiccartoons}.

\begin{figure}[ht!]
\centering

\subfloat[Barbara image]{%
\includegraphics[height=4.5cm,width=4.7cm]{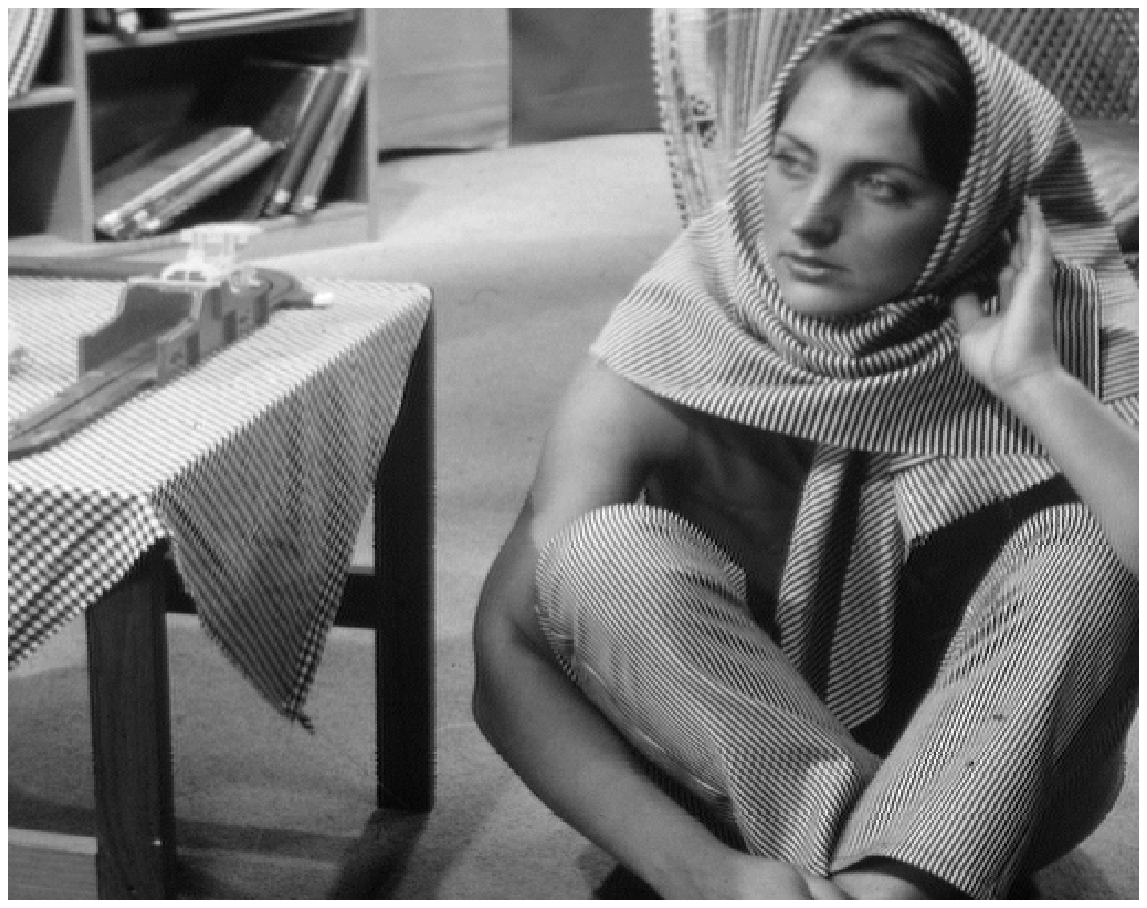}
\label{fig:cartoontexture}}

\subfloat[Cartoon Component]{%
\includegraphics[height=4.5cm,width=4.7cm]{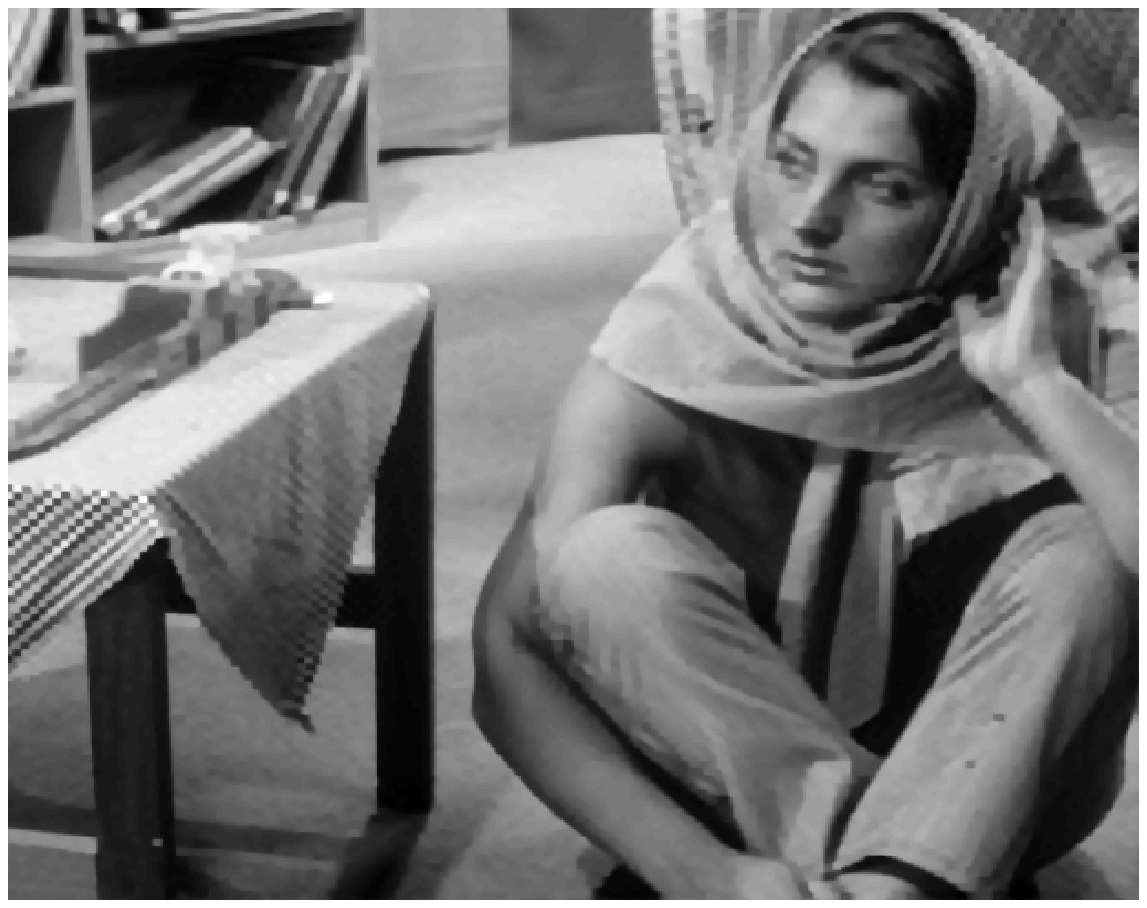}
\label{fig:cartoon} } \quad
\subfloat[Texture Component]{%
\includegraphics[height=4.5cm,width=4.7cm]{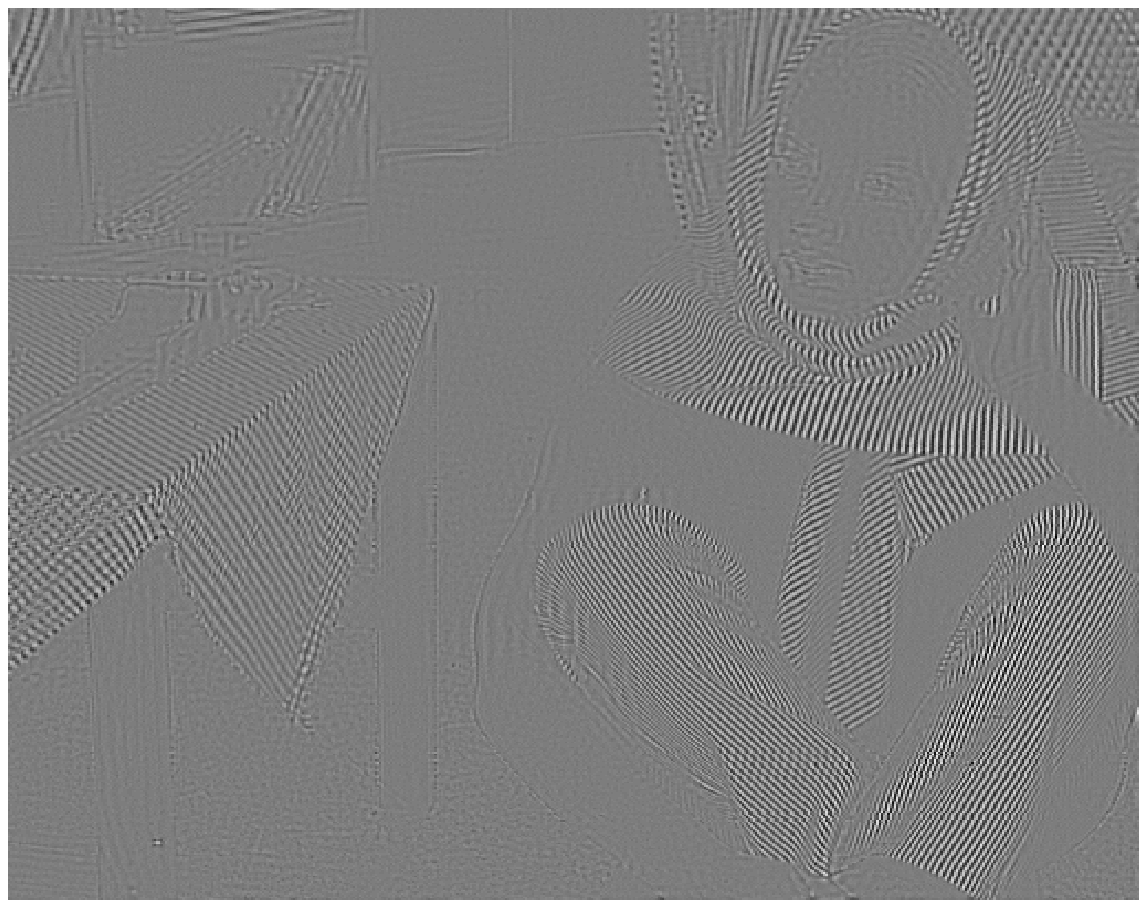}
\label{fig:texture} }


\vspace*{0.25cm}

\caption{Separation of the Barbara image into cartoon and texture
using MCALab.} \label{fig:barbara}
\end{figure}

As sparsifying systems, again curvelets or shearlets are suitable
for the cartoon part, whereas discrete cosines or a Gabor system can
be used for the texture part. MCALab uses for this separation task a
dictionary composed of curvelets and discrete cosines, see
\cite{SED05}. For illustrative purposes, we display in Figure
\ref{fig:barbara} the separation of the Barbara image into cartoon
and texture component performed by MCALab. As can be seen, all
periodic structure is captured in the texture part, leaving the
remainder to the cartoon component.

\subsection{Theoretical Results}
\label{subsec:theoryimages}

The first theoretical result explaining the successful empirical
performance of Morphological Component Analysis was derived in
\cite{DK10} by considering the separation of point- and curvelike
features in images coined the {\em Geometric Separation Problem}.
The analysis in this paper has three interesting features. Firstly,
it introduces the notion of cluster coherence (cf. Definition
\ref{def:cluster}) as a measure for the geometric arrangements of
the significant coefficients and hence the encoding of the
morphological difference of the components. It also initiates the
study of $\ell_1$ minimization in frame settings, in particular
those where singleton coherence within one frame may be high.
Secondly, it provides the first analysis of a continuum model in
contrast to the previously studied discrete models which obscure
continuum elements of geometry. And thirdly, it explores microlocal
analysis to understand heuristically why separation might be
possible and to organize a rigorous analysis. This general approach
applies in particular to two variants of geometric separation
algorithms. One is based on tight frames of radial wavelets and
curvelets and the other uses orthonormal wavelets and shearlets.

These results are today the only results providing a theoretical
foundation to image separation using ideas from sparsity
methodologies. The same situation -- separating point- and curvelike
objects --  is also considered in \cite{DK10b} however using
thresholding as a separation technique. Finally, we wish to mention
that some initial theoretical results on the separation of cartoon
and texture in images are contained in \cite{DK10c}.

Let us now dive into the analysis of \cite{DK10}. As a mathematical
model for a composition of point- and curvelike structures, the
following two components are considered: The function $\cP$ on
$\mathbb{R}^2$, which is smooth except for point singularities and
defined by
\[
\cP = \sum_{i=1}^{P}|x-x_i|^{-3/2},
\]
serves as a model for the pointlike objects, and  the distribution
$\cC$ with singularity along a closed curve $\tau:[0,1] \rightarrow
\mathbb{R}^2$ defined by
\[
\cC = \int \delta_{\tau(t)}dt,
\]
models the curvelike objects. The general model for the considered
situation is then the sum of both, i.e.,
\begin{equation}\label{eq:point_curve}
f = \cP+\cC,
\end{equation}
and the {\em Geometric Separation Problem} consists of recovering
$\cP$ and $\cC$ from the observed signal $f$.

\begin{figure}[ht!]
\centering

\mbox{\subfloat[Radial wavelets]{%
\includegraphics[scale=0.4]{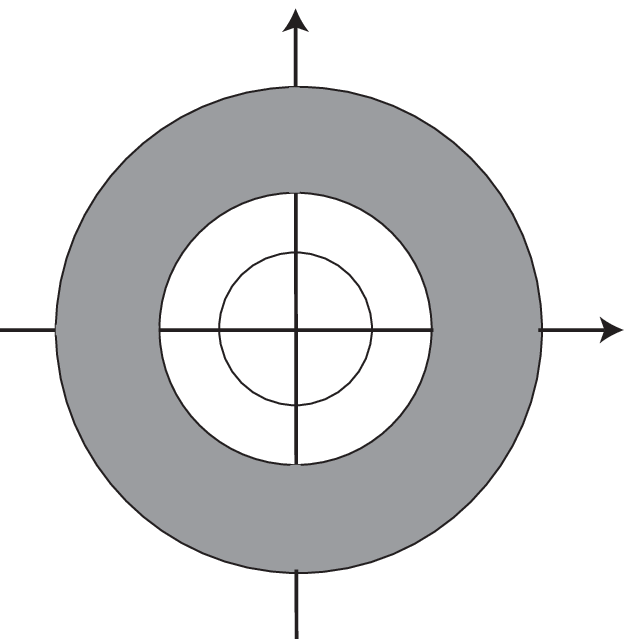}
\label{fig:radialwavelets}}

\hspace*{2.5cm}

\subfloat[Curvelets]{%
\includegraphics[scale=0.4]{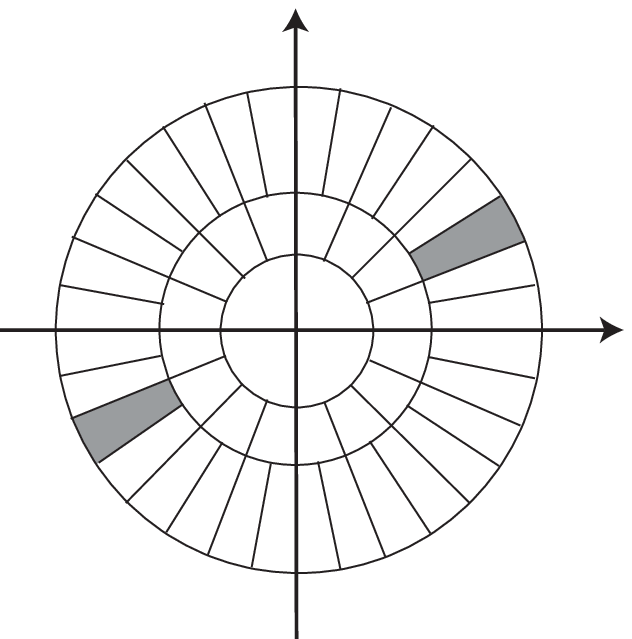}
\label{fig:curvelets1}} }

\vspace*{0.25cm}

\caption{Tiling of the frequency domain by radial wavelets and
curvelets.} \label{fig:waveletscurvelets}
\end{figure}

As discussed before, one possibility is to set up the minimization
problem using an overcomplete system composed of wavelets and
curvelets. For the analysis, radial wavelets are used due to the
fact that they provide the same subbands as curvelets. To be more
precise, let $W$ be an appropriate window function. Then {\em radial
wavelets} at scale $j$ and spatial position $k = (k_1,k_2)$ are
defined by the Fourier transforms
\[
      \hat{\psi}_{\lambda}(\xi) =  2^{-j} \cdot W(|\xi|/2^{j}) \cdot e^{ i \ip{k}{\xi /2^j}},
\]
where $\lambda = (j,k)$ indexes scale and position. For the same
window function $W$ and a `bump function' $V$, {\em curvelets} at
scale $j$, orientation $\ell$, and spatial position $k = (k_1,k_2)$
are defined by the Fourier transforms
\[
      \hat{\gamma}_{\eta}(\xi) =  2^{-j\frac{3}{4}}  \cdot W(|\xi|/2^{j}) V((\omega-\theta_{j,\ell})2^{j/2})
         \cdot e^{ i (R_{\theta_{j,\ell}}A_{2^{-j}}k)' \xi },
\]
where $\theta_{j,\ell} = 2\pi \ell /2^{j/2}$, $R_\theta$ is planar
rotation by $- \theta$ radians, $A_a$ is anisotropic scaling with
diagonal $(a, \sqrt{a})$, and we let $\eta = (j,\ell,k)$ index
scale, orientation, and scale; see \cite{CD05} for more details. The
tiling of the frequency domain generated by these two systems is
illustrated in Figure \ref{fig:waveletscurvelets}.

By using again the window $W$, we define the family of filters $F_j$
by their transfer functions
\[
    \hat{F}_j(\xi) = W(|\xi|/2^{j}), \qquad \xi \in \mathbb{R}^2 .
\]
These filters provide a decomposition of any distribution $g$ into
pieces $g_j$ with different scales, the piece $g_j$ at subband $j$
generated by filtering $g$ using $F_j$:
\[
   g_j = F_j \star g.
\]
A proper choice of $W$ then enables reconstruction of $g$ from these
pieces using the formula
\[
   g = \sum_j F_j \star g_j.
\]
Application of this filtering procedure to the model image $f$ from
\eqref{eq:point_curve} yields the decompositions
\[
f_j = F_j \star f = F_j \star (\cP + \cC) = \cP_j + \cC_j,
\]
where $(f_j)_j$ is known, and we aim to extract $(\cP_j)_j$ and
$(\cC_j)_j$. We should mention at this point that, in fact, the pair
$(\cP, \cC)$ was chosen in such a way that $\cP_j$ and $\cC_j$ have
the same energy for each $j$, thereby making the components
comparable as we go to finer scales and the separation challenging
at {\em each} scale.

Let now $\Phi_1$ and $\Phi_2$ be the tight frame of radial wavelets
and curvelets, respectively. Then, for each scale $j$, we consider
the $\ell_1$ minimization problem $(\mbox{Sep}_{a})$ stated in
\eqref{eq:newmingeneral}, which now reads:
\begin{equation}\label{eq:sep1}
\min_{P_j,C_j}\|\Phi^T_1 P_j\|_1+\|\Phi^T_2 C_j\|_1 \quad
\textrm{s.t.} \quad  f_j = P_j + C_j.
\end{equation}
Notice that we use the `analysis version' of the minimization
problem, since both radial wavelets as well as curvelets are
overcomplete systems.

The theoretical result of the precision of separation of $f_j$ via
(\ref{eq:sep1}) proved in \cite{DK10} can now be stated in the
following way:

\begin{theorem}[\cite{DK10}]
\label{theo:geo_sep} Let $\hat P_j$ and $\hat C_j$ be the solutions
to the optimization problem (\ref{eq:sep1}) for each scale $j$. Then
we have
\[
\frac{\|\cP_j - \hat P_j\|_2+\|\cC_j - \hat
C_j\|_2}{\|\cP_j\|_2+\|\cC_j\|_2} \rightarrow 0, \quad j \rightarrow
\infty.
\]
\end{theorem}

This result shows that the components $\cP_j$ and $\cC_j$ are
recovered with asymptotically arbitrarily high precision at very
fine scales. The energy in the pointlike component is completely
captured by the wavelet coefficients, and the curvelike component is
completely contained in the curvelet coefficients. Thus, the theory
evidences that the Geometric Separation Problem can be
satisfactorily solved by using a combined dictionary of wavelets and
curvelets and an appropriate $\ell_1$ minimization problem, as
already the empirical results indicate.

We next provide a sketch of proof and refer to \cite{DK10} for the
complete proof.\\[-0.3cm]

\noindent {\bf Proof [Sketch of proof of Theorem \ref{theo:geo_sep}].}
The main goal will be to apply Theorem \ref{theo:cluster2} to each
scale and prove that the sequence of bounds
$\frac{2\delta}{1-2\mu_c}$ converges to zero. For this, let $j$ be
arbitrarily fixed, and apply Theorem \ref{theo:cluster2} in the
following way: \bitem
\item $S$: Filtered signal $f_j$ $( = \cP_j + \cC_j)$.
\item $\Phi_1$: Wavelets filtered with $F_j$.
\item $\Phi_2$: Curvelets filtered with $F_j$.
\item $\Lambda_1$: Significant wavelet coefficients of $\cP_j$.
\item $\Lambda_2$: Significant curvelet coefficients of $\cC_j$.
\item $\delta_j$: Degree of approximation by significant coefficients.
\item $(\mu_c)_j$: Cluster coherence of wavelets-curvelets.
\eitem If \beq \label{eq:geo_sep1} \frac{2\delta_j}{1-2(\mu_c)_j} =
o(\|\cP_j\|_2+\|\cC_j\|_2)\quad \mbox{as } j \to \infty \eeq can be
then shown, the theorem is proved.

One main problem to overcome is the highly delicate choice of
$\Lambda_1$ and $\Lambda_2$. It would be ideal to define those sets
in such a way that \beq \label{eq:geo_sep2} \delta_j =
o(\|\cP_j\|_2+\|\cC_j\|_2)\quad \mbox{as } j \to \infty \eeq and
\beq \label{eq:geo_sep3} (\mu_c)_j \to 0 \quad \mbox{as } j \to
\infty \eeq are true. This would then imply \eqref{eq:geo_sep1},
hence finish the proof.

A microlocal analysis viewpoint now provides insight into how to
suitably choose $\Lambda_1$ and $\Lambda_2$ by considering the
wavefront sets of $\cP$ and $\cC$ in phase space $\mathbb{R}^2
\times [0,2\pi)$, i.e.,
\[
      WF(\cP) = \{x_i\}_{i=1}^P \times [0,2\pi)
\]
and
\[
      WF(\cC) =  \{ (\tau(t), \theta(t)): t \in [0,L(\tau)] \},
\]
where $\tau(t)$ is a unit-speed parametrization of $\tau$ and
$\theta(t)$ is the normal direction to $\tau$ at $\tau(t)$.
Heuristically, the significant wavelet coefficients should be
associated with wavelets whose index set is `close' to $WF(\cP)$ in
phase space and, similarly, the significant curvelet coefficients
should be associated with curvelets whose index set is `close' to
$WF(\cC)$. Thus, using Hart Smith's phase space metric,
\begin{eqnarray*}
d_{HS}((b,\theta);(b',\theta')) & = & \absip{e_\theta}{b-b'} +
\absip{e_{\theta'}}{b-b'} + |b-b'|^2 + |\theta-\theta'|^2,
\end{eqnarray*}
where $e_\theta = (\cos(\theta), \sin(\theta))$, an `approximate'
form of sets of significant wavelet coefficients is
\[
\Lambda_{1,j} = \{\mbox{wavelet lattice}\} \cap \{(b,\theta) :
d_{HS}((b,\theta);WF(\cP)) \le \eta_j 2^{-j}\},
\]
and an `approximate' form of sets of significant curvelet
coefficients is
\[
\Lambda_{2,j} = \{\mbox{curvelet lattice}\} \cap \{(b,\theta) :
d_{HS}((b,\theta);WF(\cC)) \le \eta_j 2^{-j}\}
\]
with a suitable choice of the distance parameters $(\eta_j)_j$. In
the proof of Theorem \ref{theo:geo_sep}, the definition of
$(\Lambda_{1,j})_j$ and $(\Lambda_{2,j})_j$ is much more delicate,
but follows this intuition. Lengthy and technical estimates then
lead to \eqref{eq:geo_sep2} and \eqref{eq:geo_sep3}, which -- as
mentioned before -- completes the proof.
\qed

Since it was already mentioned in Subsection
\ref{subsec:empiricalimages} that a combined dictionary of wavelets
and shearlets might be preferable, the reader will wonder whether
the just discussed theoretical results can be transferred to this
setting. In fact, this is proven in \cite{Kut10}, see also
\cite{DK09}. It should be mentioned that one further advantage of
this setting is the fact that now a basis of wavelets can be
utilized in contrast to the tight frame of radial wavelets explored
before.

\begin{figure}[ht!]
\centering

\mbox{\subfloat[Wavelets]{%
\includegraphics[scale=0.4]{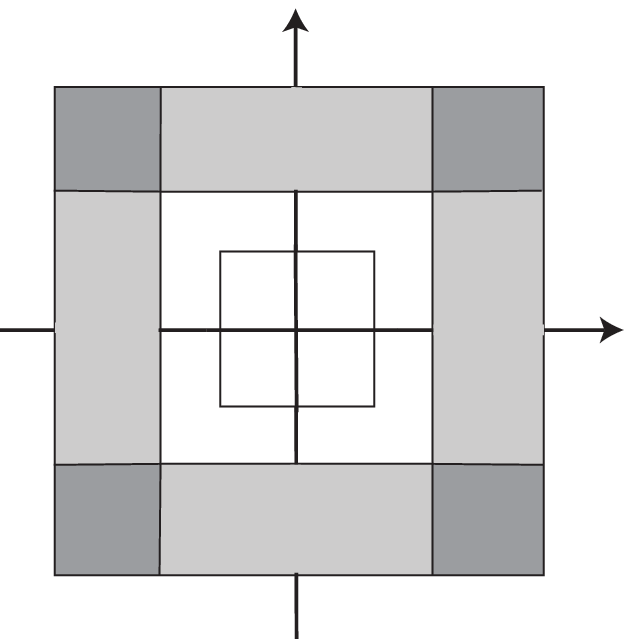}
\label{fig:ONwavelets}}

\hspace*{2.5cm}

\subfloat[Shearlets]{%
\includegraphics[scale=0.4]{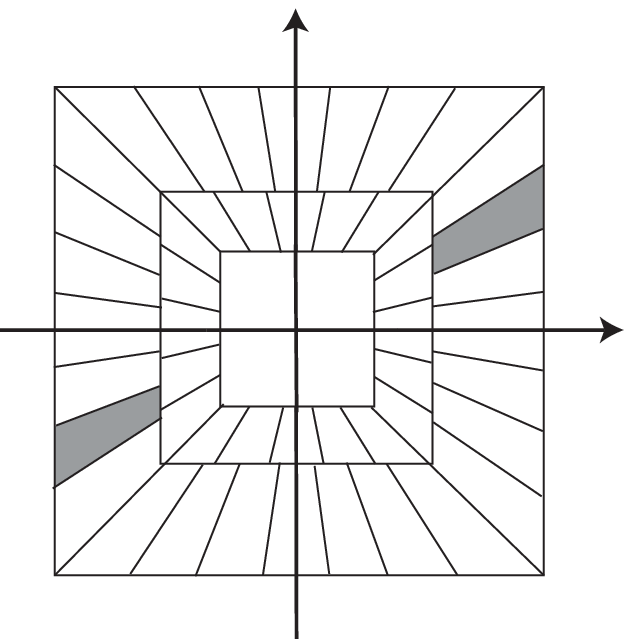}
\label{fig:shearlets1}} }

\vspace*{0.25cm}

\caption{Tiling of the frequency domain by orthonormal Meyer
wavelets and shearlets.} \label{fig:waveletsshearlets}
\end{figure}

As a wavelet basis, we now choose {\em orthonormal Meyer wavelets},
and refer to \cite{Mal98} for the definition. For the definition of
shearlets, for $j \ge 0$ and $k \in \mathbb{Z}$, let --~the notion
$A_{2^j}$ was already introduced in the definition of curvelets --
$\tilde{A}_{2^j}$ and $S_k$ be defined by
\[
\tilde{A}_{2^j} = \begin{pmatrix} 2^{j/2} & 0 \\ 0 & 2^{j}
\end{pmatrix} \quad \text{and} \quad S_k = \begin{pmatrix} 1 & k \\
0 & 1 \end{pmatrix}.
\]
For $\phi, \psi, \tilde{\psi} \in L^2({\mathbb R}^2)$, the {\em
cone-adapted discrete shearlet system} is then the union of
\[
\{\phi(\cdot-m) : m \in {\mathbb Z}^2\},
\]
\[
\{2^{\frac34 j} \psi(S_kA_{2^j}\, \cdot \, -m) : j \ge 0, -\lceil
2^{j/2}\rceil \leq k \leq \lceil 2^{j/2}\rceil, m \in {\mathbb
Z}^2\},
\]
and
\[
\{2^{\frac34 j} \tilde{\psi}(S^T_k\tilde{A}_{2^j}\, \cdot \, -m) : j
\ge 0, -\lceil 2^{j/2}\rceil \leq k \leq \lceil 2^{j/2}\rceil, m \in
{\mathbb Z}^2 \}.
\]
The term `cone-adapted' originates from the fact that these systems
tile the frequency domain in a cone-like fashion; see Figure
\ref{fig:shearlets1}.

As can be seen from Figure \ref{fig:waveletsshearlets}, the subbands
associated with orthonormal Meyer wavelets and shearlets are the
same. Hence a similar filtering into scaling subbands can be
performed as for radial wavelets and curvelets.

Adapting the optimization problem (\ref{eq:sep1}) by using wavelets
and shearlets instead of radial wavelets and curvelets generates
purported point- and curvelike objects $\hat W_j$ and $\hat S_j$,
say, for each scale $j$. Then the following result, which shows
similarly successful separation as Theorem \ref{theo:geo_sep}, was
derived in \cite{Kut10} with the new concept of sparsity
equivalence, here between shearlets and curvelets, introduced in the
same paper as main ingredient.

\begin{theorem}[\cite{Kut10}]
We have
\[
\frac{\|\cP_j - \hat W_j\|_2+\|\cC_j - \hat
S_j\|_2}{\|\cP_j\|_2+\|\cC_j\|_2} \rightarrow 0, \quad j \rightarrow
\infty.
\]
\end{theorem}






\begin{thebibliography}{99}

\bibitem{AEB06}
Aharon, M., Elad, M., and Bruckstein, A.M. (2006).
\newblock The K-SVD: An algorithm for designing of overcomplete dictionaries for sparse representation,
{\em IEEE Trans. Signal Proc.}, {\bf 54(11)}, 4311--4322.

\bibitem{BSFMD07}
Bobin, J., Starck, J.-L., Fadili, M.J., Moudden, Y., and Donoho,
D.L. (2007).
\newblock Morphological component analysis: An adaptive thresholding strategy,
{\em IEEE Trans. Image Proc.}, {\bf 16(11)}, 2675--2681.

\bibitem{BGN08}
Borup, L., Gribonval, R., and  Nielsen, M. (2008).
\newblock Beyond coherence: Recovering structured time-frequency representations,
{\em Appl. Comput. Harmon. Anal.}, {\bf 24(1)}, 120--128.

\bibitem{BDE09}
 Bruckstein, A.M., Donoho, D.L., and  Elad, M. (2009).
\newblock From sparse solutions of systems of equations to sparse modeling of signals and images,
SIAM Review, {\bf 51(1)}, 34--81.

\bibitem{CD05}
Cand\`es, E.J. and Donoho, D.L. (2005).
\newblock Continuous curvelet transform: II. Discretization of frames,
{\em Appl. Comput. Harmon. Anal.}, {\bf 19(2)}, 198--222.

\bibitem{CDS98}
Chen, S. S., Donoho, D. L., and Saunders, M. A. (1998).
\newblock Atomic decomposition by basis pursuit,
{\em SIAM J. Sci. Comput.}, {\bf  20(1)}, 33--61.

\bibitem{CW93}
 Coifman, R.R. and  Wickerhauser, M.V. (1993).
\newblock Wavelets and adapted waveform analysis. A toolkit for signal processing and numerical analysis,
{\em Different perspectives on wavelets} (San Antonio, TX, 1993),
119--153, Proc. Sympos. Appl. Math., {\bf 47}, Amer. Math. Soc.,
Providence, RI.


\bibitem{Don01}
Donoho, D.L. (2001).
\newblock Sparse components of images and optimal atomic decomposition,
{\em Constr. Approx.}, {\bf 17(3)}, 353--382.

\bibitem{Don06c}
 Donoho, D.L. (2006).
\newblock Compressed sensing,
{\em IEEE Trans. Inform. Theory}, {\bf 52(4)},  1289--1306.

\bibitem{DE03}
Donoho, D.L. and Elad, M. (2003).
\newblock Optimally sparse representation in general (non\-or\-tho\-gonal) dictionaries via $l\sp 1$ minimization,
{\em Proc. Natl. Acad. Sci. USA}, {\bf 100(5)},  2197--2202.

\bibitem{DH01}
Donoho, D.L.  and Huo, X. (2001).
\newblock Uncertainty principles and ideal atomic decomposition,
{\em IEEE Trans. Inform. Theory}, {\bf 47(7)}, 2845--2862.

\bibitem{DK09}
Donoho, D.L. and Kutyniok, G. (2009).
\newblock Geometric separation using a wavelet-shearlet dictionary,
{\em SampTA'09} (Marseille, France, 2009), Proc., 2009.

\bibitem{DK10b}
Donoho, D.L. and Kutyniok, G. (2010).
\newblock Geometric separation by single-pass alternating thresholding,
preprint.

\bibitem{DK10}
Donoho, D.L. and Kutyniok, G. (2010).
\newblock Microlocal analysis of the geometric separation problem,
preprint.

\bibitem{DK10c}
Donoho, D.L. and Kutyniok, G. (2011).
\newblock Geometric separation of cartoons and texture via $\ell_1$ minimization,
preprint.

\bibitem{DS89}
 Donoho, D.L. and  Stark, P.B. (1989).
\newblock Uncertainty principles and signal recovery,
{\em SIAM J. Appl. Math.}, {\bf 49(3)}, 906--931.

\bibitem{DS09}
Duarte-Carvajalino, J.M. and Sapiro, G. (2009).
\newblock Learning to sense sparse signals: Simultaneous sensing matrix and
sparsifying dictionary optimization, {\em IEEE Trans. Image Proc.},
{\bf 18(7)}, 1395--1408.

\bibitem{Ela10}
Elad, M. (2010).
\newblock {\em Sparse and redundant representations},
Springer, New York.

\bibitem{EB02}
Elad, M. and Bruckstein, A. M. (2002).
\newblock A generalized uncertainty principle and sparse representation in pairs of bases,
{\em IEEE Trans. Inform. Theory}, {\bf 48(9)}, 2558--2567.

\bibitem{ESQD05}
Elad, M.,  Starck, J.-L., Querre, P., and  Donoho, D.L. (2005).
\newblock Simultaneous cartoon and texture image inpainting using morphological component analysis (MCA),
{\em Appl. Comput. Harmon. Anal.},  {\bf 19(3)},  340--358.

\bibitem{EAHH99}
Engan, K., Aase, S.O., and Hakon-Husoy, J.H. (1999).
\newblock Method of optimal directions for frame design,
{\em IEEE Int. Conf. Acoust., Speech, Signal Process.}, {\bf 5},
2443-–2446.

\bibitem{GB03}
 Gribonval, R. and  Bacry, E. (2003).
\newblock Harmonic decomposition of audio signals with matching pursuit,
{\em IEEE Trans. Signal Proc.}, {\bf 51(1)}, 101--111.

\bibitem{GN03}
 Gribonval, R. and  Nielsen, M. (2003).
\newblock Sparse representations in unions of bases,
{\em IEEE Trans. Inform. Theory},  {\bf 49(12)},  3320--3325.

\bibitem{GKL06}
 Guo, K.,  Kutyniok, G., and  Labate, D. (2006).
\newblock Sparse multidimensional representations using anisotropic dilation and shear operators,
{\em Wavelets and Splines} (Athens, GA, 2005), Nashboro Press,
Nashville, TN, 2006, 189--201.

\bibitem{KT09}
Kowalski, M.  and  Torr\'{e}sani, B. (2010).
\newblock Sparsity and persistence: Mixed norms provide simple signal models with dependent coefficients,
{\em Signal, Image and Video Proc.}, to appear.

\bibitem{Kut10}
 Kutyniok, G. (2010).
\newblock Sparsity equivalence of anisotropic decompositions,
preprint.

\bibitem{KLL10}
 Kutyniok, G.,  Lemvig, J., and Lim, W.-Q (2010).
\newblock Compactly supported shearlets,
{\em Approximation Theory XIII} (San Antonio, TX, 2010), Springer,
to appear.

\bibitem{KL10a}
Kutyniok, G. and  Lim, W.-Q (2010).
\newblock Compactly supported shearlets are optimally sparse,
preprint.

\bibitem{KL10b}
Kutyniok, G. and  Lim, W.-Q (2010).
\newblock Image separation using shearlets,
preprint.

\bibitem{Mal98}
Mallat, S.G. (1998).
\newblock {\em A wavelet tour of signal processing},
Academic Press, Inc., San Diego, CA.

\bibitem{MZ93}
 Mallat, S.G. and Zhang,  Z. (1993).
\newblock Matching pursuits with time-frequency dictionaries,
{\em IEEE Trans. Signal Proc.}, {\bf 41(12)}, 3397--3415.

\bibitem{MAC02}
  Meyer, F.G., Averbuch, A., and  Coifman, R.R. (2002).
\newblock Multi-layered image representation: Application to image compression,
 {\em IEEE Trans. Image Proc.}, {\bf 11(9)}, 1072--1080.

\bibitem{Sta10}
Starck, J.-L., Murtagh, F., and Fadili, J.M. (2010).
\newblock {\em Sparse Image and Signal Processing: Wavelets, Curvelets, Morphological Diversity},
Cambridge University Press, New York, NY.

\bibitem{SED04}
Starck, J.-L., Elad, M., and  Donoho, D.L. (2005).
\newblock Redundant multiscale transforms and their application for morphological component analysis,
{\em Adv. Imag. Electr. Phys.}, {\bf 132}, 287--348.

\bibitem{SED05}
Starck, J.-L., Elad, M., and  Donoho, D.L. (2005).
\newblock Image decomposition via the combination of sparse representations and a variational approach,
{\em IEEE Trans. Image Proc.}, {\bf 14(10)}, 1570--1582.

\bibitem{SMBED05}
 Starck, J.-L., Moudden, Y., Bobin, J., Elad, M., and  Donoho, D.L. (2005).
\newblock Morphological component analysis,
{\em Wavelets XI} (San Diego, CA, 2005), SPIE Proc. 5914, SPIE,
Bellingham, WA.

\bibitem{Tro04}
Tropp, J.A. (2004).
\newblock  Greed is good: Algorithmic results for sparse approximation,
{\em IEEE Trans. Inform. Theory},  {\bf 50(10)},  2231--2242.

\bibitem{Tro08}
Tropp, J.A. (2008).
\newblock  On the linear independence of spikes and sines,
{\em J. Fourier Anal. Appl.},  {\bf 14(5-6)},  838--858.

\bibitem{Tro10}
Tropp, J.A. (2010).
\newblock  The sparsity gap: Uncertainty principles proportional to dimension,
{\em Proc. 44th IEEE Conf. Information Sciences and Systems} (CISS),
1--6, Princeton, NJ.

\end{thebibliography}
\end{document}